\documentstyle[11pt]{article}
\normalbaselines
\begin{document}
\bibliographystyle{unsrt}

\def\bea*{\begin{eqnarray*}}
\def\eea*{\end{eqnarray*}}
\def\ba{\begin{array}}
\def\ea{\end{array}}
% --------------------------------------------------------------
% If you want to see the names of the equations and references
% set below \count1=0 otherwise \count1=1
% --------------------------------------------------------------
\count1=1
% --------------------------------------------------------------
\def\be{\ifnum \count1=0 $$ \else \begin{equation}\fi}
\def\ee{\ifnum\count1=0 $$ \else \end{equation}\fi}
\def\ele(#1){\ifnum\count1=0 \eqno({\bf #1}) $$ \else \label{#1}\end{equation}\fi}
\def\req(#1){\ifnum\count1=0 {\bf #1}\else \ref{#1}\fi}
\def\bea(#1){\ifnum \count1=0   $$ \begin{array}{#1}
\else \begin{equation} \begin{array}{#1} \fi}
\def\eea{\ifnum \count1=0 \end{array} $$
\else  \end{array}\end{equation}\fi}
\def\elea(#1){\ifnum \count1=0 \end{array}\label{#1}\eqno({\bf #1}) $$
\else\end{array}\label{#1}\end{equation}\fi}
\def\cit(#1){
\ifnum\count1=0 {\bf #1} \cite{#1} \else 
\cite{#1}\fi}
\def\bibit(#1){\ifnum\count1=0 \bibitem{#1} [#1    ] \else \bibitem{#1}\fi}
\def\ds{\displaystyle}
\def\hb{\hfill\break}
\def\comment#1{\hb {***** {\em #1} *****}\hb }

\newcommand{\TZ}{\hbox{T\hspace{-5pt}T}}
\newcommand{\MZ}{\hbox{I\hspace{-2pt}M}}
\newcommand{\ZZ}{\hbox{Z\hspace{-3pt}Z}}
\newcommand{\NZ}{\hbox{I\hspace{-2pt}N}}
\newcommand{\RZ}{\hbox{I\hspace{-2pt}R}}
\newcommand{\CZ}{\,\hbox{I\hspace{-6pt}C}}
\newcommand{\PZ}{\hbox{I\hspace{-2pt}P}}
\newcommand{\QZ}{\hbox{I\hspace{-6pt}Q}}
\newcommand{\KZ}{\hbox{I\hspace{-2pt}K}}

\vbox{\vspace{38mm}}
\begin{center}
{\LARGE \bf  Picard Groups of Hypersurfaces in Toric Varieties }\\[5mm]

Shi-shyr Roan\footnote{Supported in part by the NSC grant of Taiwan.}\\{\it Institute of Mathematics \\ Academia Sinica \\ 
Taipei , Taiwan \\ (e-mail: maroan@ccvax.sinica.edu.tw)} \\[5mm]
\end{center}

\begin{abstract} We study the structure of rational 
Picard groups of hypersurfaces of toric varieties. By using 
the fan structure associated to the ambient toric variety,
an explicit basis of the Picard group is described by certain
combinatorial data. We shall also discuss the application to Calabi-Yau
spaces.  
\end{abstract}
\par \vspace{5mm} \noindent

\section{Introduction}
Throughout this paper by an orbifold we shall
mean a complex variety with at worst abelian quotient 
singularities. The orbifolds we shall study here are
toric varieties and their quasi-smooth hypersurfaces. 
Much of what was needed in toric geometry is already in 
\cite{Dan} \cite{F} \cite{KKMS} \cite{O78} \cite{O88}, 
and the present paper borrows heavily from those 
references. 
The purpose of this paper is to determine the (rational) Picard
groups, i.e. the groups of $\QZ$-divisors, of hypersurfaces in
toric varieties. The
quantitative aspect of Picard groups, as well as the 
qualitative one,  will be the main concern of this work. In 
algebraic geometry, the toric divisors  provide
a combinatorial description for the structure of 
$\QZ$-divisors of a toric variety \cite{BC} \cite{MT}. However, one
can also study this subject from sympletic geometry point of
view \cite{Au} \cite{Gui}. As in the case of projective spaces 
where line bundles are obtained by  
the induced bundles from Hopf bundle through characters of ${\CZ}^*$, 
we can also describe the Picard groups of toric varieties in a similar
manner. The present work focuses on the structure
on Picard groups of hypersurfaces in toric varieties,
which ultimately yields some combinatorial basis-representation
for $H^{1, 1}$ of Calabi-Yau (CY) mirror 3-folds.
There has been considerable progress on this subject over the last few 
years \cite{Y}: the special case on quotients of 
"Fermat-type" CY hypersurfaces in \cite{R91} 
\cite{R93y}; and also the general cases but focusing on Picard
numbers only in \cite{B}. The goal of this sequel is to show that
an investigation of the general properties reveals a remarkably 
profound structure, and hence to enrich our understanding on toric 
geometry and mirror symmetry.

The following is a summary of the contents of this paper:
Sect. 2 contains some basic definitions in this work, 
including orbifold principal bundles
with torus as the structure group. The homogeneous coordinate 
system of toric varieties found by Audin \cite{Au} and Cox
\cite{Cox} \cite{O94} \cite{AGM} is the main example for such notion in
the present paper. In Sect. 3, the  Picard groups of 
toric varieties are constructed through
the induced bundle-construction which behaves much the same 
as in the ususal projective space. In Sect. 4, we aim to give
an overview of some basic properties on quasi-smooth hypersurfaces
in toric varieties for later use.
For the purpose of illustrating the general results in the 
previous sections, we shall present some standard
examples of toric varieties in Sect. 5 by connecting them 
with some well-known facts on those varieties.
In the next two sections we shall deal with the most 
significant objects in this paper. In Sect. 6, Picard groups
of hypersurfaces of toric varieties are determined through
the related combinatorial data. The hypersurfaces are in 
general not ample, however we must assume some additional
(somewhat unpleasant) conditions on the toric varieties.
The prime examples in our mind are CY
spaces appeared in the context of mirror symmetry
\cite{B} \cite{R91} \cite{R93y} \cite{Y}, whose structure
has led us to this work. In Sect. 7, we shall apply the results
of previous sections to anti-canonical hypersurfaces
in the toric varieties of reflexive polytope-type introduced
by Batyrev in \cite{B}. With an additional "maximal" condition
on the combinatorial data, the hypersurfaces are smooth in
the dimension 3. Hence we obtain a combinatorial 
descrpition of Picard groups for these CY spaces, which 
generalize the results in \cite{R91} \cite{R93y}
on the quotients of Fermat-type CY hypersurfaces.
In this paper, we have endeavored to put the presentation in the context
of complex geometry in the hope 
of making easier access to differential geometers and topologists.
Some of what is included should be undoubtedly known to experts in toric
geometry; nonetheless, it seems difficult to find an appropriate
references, especially in Sect. 2 and 4, some subjects we have
taken pains to include. We have put some elementary arguments 
for a few well-known facts of the content in the appendix for 
easy reference.

\par \vspace{0.2in} \noindent
{\bf Notations } 
\par \vspace{0.2in} \noindent 
We prepare some notations for easier presentation of 
this paper. \par  \noindent

$ A_{\KZ} = A \bigotimes_{\ZZ} {\KZ} \ \ , \ \ \mbox{for
\ an \ abelian \ group \ } A \mbox{ \ and \ a field \ }\ 
{\KZ} = {\QZ} , {\RZ} , {\CZ} \ . $
\par \noindent
For a $n-$dimensional lattice $L$, we denote

${\TZ} ( L ) =  L_{\CZ} / L \ \
\ ( \ = \mbox{exp} ( 2 \pi i L_{\CZ} ) \ ) \ , $ \ \ the $n-$dimensional 
(algebraic) torus having $L$ as the group of one-parameter subgroups.

$L^* = \ \mbox{Hom} ( L , \ZZ ) $ , \ the dual lattice of $L$ .

$< * , * > \ : \ L_{\RZ} \times L^*_{\RZ} \longrightarrow \RZ 
$ , \ the non-degenerate 
natural pairing which takes integral values on $L \times L^*$. For a
subset $\tau$ of $L_{\RZ}$, the annihilator of $\tau$ is defined
by
\[
\tau^\bot : = \{ \ y \in L^*_{\RZ} \ | \ < x , y > = 0 \ \ 
\forall \ \ x \in \tau \  \} .
\]

$\Sigma$ = \ a fan in $L_{\RZ}$ , which will always be 
a rational simplical fan for the lattice $L$ throughout the paper unless otherwise 
specified. $\Sigma^{(i)}$ denotes the $i$-th skeleton of $\Sigma$.

${\PZ}_{ ( \Sigma , L )}$ = \ the ${\TZ} ( L )$-toric variety associated
to a complete fan $\Sigma$, 
( note that ${\PZ}_{ ( \Sigma , L )}$ is an orbifold by the 
assumption on $\Sigma$ ).

\section{Orbifold principal bundles}
For convenience of later discussions, we 
introduce the following notion on orbifolds
as a generalization of principle bundles and induced line 
bundles in the class of complex manifolds:
\par \vspace{0.2in} \noindent 
{\bf Definition 1 . } Let $\TZ$ be a (algebraic) torus. 
A morphism of orbifolds,
\[
\pi : P \longrightarrow X \ ,
\]
is called a principal orbifold  $\TZ$-bundle if the following
conditions hold:

(i) $P$ is a $\TZ$-space with the right $\TZ$-action:  
\[
P \times {\TZ} \longrightarrow P ,
\]
and $\pi$ is a $\TZ$-equivariant morphism with the trivial $\TZ$-structure
on $X$. 

(ii) For $x \in X$, there is an open neighborhood
$U_x $ of $x$, and a finite abelian group $G_x$ such that  
\[
\pi_{ \mbox{rest}} \ : \ \pi^{-1} ( U_x ) \longrightarrow U_x
\]
is $\TZ$-equivalent to a left $G_x$-quotient of 
the product bundle,
\[ 
G_x \setminus ( B_{\epsilon} ( 0 ) \times 
{\TZ}   
\stackrel{pr}{\longrightarrow} B_{\epsilon} ( 0 ) ) \ , \ \
B_{\epsilon} ( 0 ) : =  
\{ z \in \CZ^n \ ; \ |z| < \epsilon \} \ , 
\]
where the $G_x$-action is given by
\[
G_x \times  B_{\epsilon} ( 0 ) \times 
{\TZ} \longrightarrow B_{\epsilon} ( 0 ) \times 
{\TZ} \ , \ ( g , ( z , t ) ) \mapsto 
g \cdot ( z , t ) = ( g \cdot z , \  \varphi (g ) t ) \ ,
\]
for some group homomorphism $\varphi : G_x \rightarrow \TZ$.

(iii) $G_x$ is the trivial group for a non-singular point 
$x$. 
$\Box$ \par \vspace{0.2in} \noindent
Note that $\TZ$ acts on $\pi^{-1} (x) $ transitively
for each $x \in X$. Hence all elements
of $\pi^{-1} (x) $ have the same finite isotropy subgroup of 
$\TZ$, denoted by $g_x$. For $U_x$ in (ii) of the above 
definition, an element $y$ of $U_x$  corresponds
to a $G_x$-orbit $[z_0]$ for some  
$z_0 \in \CZ^n$. We have the following isomorphic germs of
analytic spaces:
\[
( U_x , y ) \simeq ( \{ z \in \CZ^n \ | \ | z - z_0 | < \delta  
\} / G_{x , y } \ , \ [ z_0 ] \ ) \ ,
\]
where $G_{x , y}$ is the subgroup of $G_x$ stablizing the fiber over 
$z_0$, ( note that the group $G_{x , y}$ does not depend on the choice
of $z_0$ ). Then we have a surjective homomorphism
\[
G_{x , y} \longrightarrow g_y \ , \ \ \
\mbox{for \ } y \in U_x \ .
\]
{\bf Definition 2. } For a principal orbifold $\TZ$-bundle,
\[
\pi : P \longrightarrow X  \ ,
\]
and a character of ${\TZ}$,
\[
\chi : {\TZ} \longrightarrow {\CZ}^* \ ,
\]
$P \times_{\TZ} {\CZ}$ is the quotient  
space of $P \times {\CZ}$ by the $\TZ$-action:
\[
( p , \zeta ) \cdot t = ( p \cdot t \ , \ \chi ( t ) \zeta ) \ .
\]
The fibration
\[ 
P \times_{\TZ} {\CZ}  \longrightarrow X \ , \ \ 
[ ( p , \zeta )] \mapsto \pi ( p ) \ ,
\]
will be denoted by ${\cal O} ( \chi )$, called the 
orbifold line bundle induced by $\chi^{-1}$.
$\Box$ \par \vspace{0.2in} \noindent
Note that ${\cal O} ( \chi )$ is a line bundle 
outside the singular set of $X$. The 
sections of ${\cal O} ( \chi )$ over an open set $U$ of $X$ 
can be regarded as functions on $\pi^{-1} ( U )$. In fact, one
has
\[
\Gamma ( U , {\cal O} ( \chi ) ) = \{ f : \pi^{-1} ( U ) \rightarrow {\CZ} 
\ | \ f ( p \cdot t ) = \chi ( t ) f ( p ) \ \ \mbox{ for } \
p \in \pi^{-1} ( U ) \ , \ t \in \TZ \ \} . 
\]
We shall make this identification in what follows.
For compact $X$, ${\cal O} ( \chi )$ is considered as  
an element of the rational Picard group of $X$ by the following lemma:
\par \vspace{0.2in} \noindent  
{\bf Lemma 1.} \ For a compact variety $X$, 
${\cal O} ( \chi^N )$ is a line bundle over $X$
for some positive integer $N$.
\par \vspace{0.2in} \noindent
{\it Proof. } Let $U_x$, $G_x$ be the same as in the 
condition (ii) of Definition 1 . By the 
compactness of $X$, there is a positive interger $N$ 
divisible by $| G_x |$ for all $x \in X$. Claim : 
${\cal O} ( \chi^N )$ is a line bundle which is trivial
over $U_x$. By the construction
of orbifold line bundle, 
${\cal O} ( \chi^N )_{| U_x } $ is equivalent to
\[ G_x \setminus ( B_{\epsilon}(0)  
\times {\CZ} \ \ 
\stackrel{pr}{\longrightarrow} \ \ 
B_{\epsilon}(0) ) \ ,
\]
with the $G_x-$action given by
$$
\begin{array}{llr}
g \cdot ( z , \zeta )  &= ( g \cdot z , \ \chi^{N}(g \cdot 1) \zeta ) &\\
    &= ( g \cdot z , \  \zeta ) , &  \ \ ( 
\mbox{by} \ \ |G_x| \ | \ N \ ) \  .
\end{array}
$$
Therefore ${\cal O} ( \chi^N )_{| U_x } $ is isomorphic to $U_x \times
\CZ$. 
$\Box$ \par \vspace{0.2in} \noindent 
Let $\PZ_{( \Sigma, L)}$ be the toric variety for a 
$n$-dimensional lattice $L$ and a complete fan $\Sigma$ 
in $L_{\RZ}$. By \cite{Au} \cite{Cox}, $\PZ_{( \Sigma, L)}$ 
can be regarded
as the orbit-space  for some torus action on a "homogeneous" affine
variety. We are going to show that this fibration over
$\PZ_{( \Sigma, L)}$ is in fact a principal orbifold bundle. 
The "homogeneous" affine variety will be described through the 
toric construction as in \cite{O94}.  

Let $\Sigma^{1} $ be the collection of primitive 
elements in $L$ which are the generators of elements in 
1-skeleton $\Sigma^{(1)}$. For a simplicial cone $\sigma$ of
$\Sigma$, the simplex in $L_{\RZ}$ spanned by
$\Sigma^{1} \cap \sigma $ and the origin is denoted by 
$\triangle_{\sigma }$. Consider the (compact) polytope in 
$L_{\RZ}$,
\be
\triangle ( L ) : = \bigcup_{ \sigma \in \Sigma } 
\triangle_{\sigma } .
\ele(deftria)
The faces of $\triangle_{\sigma }$'s form a 
(simplicial ) triangulation
of $\triangle ( L )$, denoted by $\{ s_i \}_{i \in I}$, with $\Sigma^{1} \cup \{ 0 \}$ as the set of 
0-simplices. Hence  the complete fans in 
$L_{\RZ}$ are in one-to-one correspondence with the triangulated
polytopes $(\triangle ( L ), \{ s_i \}_{i \in I} )$ in 
$L_{\RZ}$ such that the following properties hold for
 0-simplices $s_i$,
$$  
\begin{array}{cll}  
s_i \in {\rm interior}(\triangle ( L )) & \Longleftrightarrow &
s_i = {\rm the \ origin \ } \vec{0} \ ,  \\
s_i \neq \vec{0} & \Longrightarrow & s_i : {\rm an 
\ primitive \ element \ in \ } L \ .
\end{array}
$$
The toric variety $\PZ_{ ( \Sigma , L )}$ is determined by 
the triangulation of $\triangle ( L )$. It is known that the 
irreducible toric divisors 
of ${\PZ}_{( \Sigma , L )}$ are in one-to-one 
correspondence with the elements in $\Sigma^{1}$, and denote 
$e_{\delta}$ the toric divisor corresponding 
to an element $\delta \in \Sigma^{1}$. Consider 
the free $\ZZ$-module  
$D_0$ of rank $d \ (: = | \Sigma^{1} |)$ 
with the $\ZZ$-basis given by  the set of all irreducible  toric divisors  of 
$ \PZ_{ ( \Sigma , L )}$, and we shall  
denote $\ZZ$-basis elements of $D_0$ by $ e_{\delta}$ again if no 
confusion could arise, 
$$
D_0 : =  \bigoplus_{ \delta \in \Sigma^{1} } 
{\ZZ} e_{\delta} \ .  
$$
( Note that no linear equivalent relation of divisors 
exists among $e_{\delta}$'s in $D_0$.)
Define the homomorphism
\[
\beta_0 : D_0 \
\longrightarrow L \ \ \ \mbox{with} \ \ \ \beta_0 ( e_{\delta} )
= \delta \ ,
\]
and denote its kernel by
\be
{\bf n}  \ ( = {\bf n} ( \Sigma , L ) ) : =  \mbox{Ker} (\beta_0) .
\ele(defn)
Since the cokernel of $\beta_0$ is finite ,
by tensoring with ${\RZ}$, we have the surjective linear 
map:
\be
\beta = ( \beta_0)_{\RZ} :  (D_0)_{\RZ} = \bigoplus_{ \delta \in \Sigma^{1} } {\RZ} e_{\delta} 
\longrightarrow L_{\RZ} \ .
\ele(defbeta)
Let $D$ be a lattice in $\bigoplus_{ \delta \in \Sigma^{1} } 
{\RZ} e_{\delta}$ with the following property:

(i) $ D \supseteq D_0 $.

(ii) $\beta ( D ) = L$, and $\beta$ induces an isomorphism of finite abelian groups:
\[
D / D_0
\simeq L / L_0 \ \ , \ \ \mbox{with} \ \ 
L_0 :=  \sum_{ \delta \in \Sigma^{1} } {\ZZ} \delta  .
\]
Note that $D$ is equal to
$\bigoplus_{ \delta \in \Sigma^{1} } {\ZZ} e_{\delta} $
if $\beta_0$ is surjective. In general,
the lattice $D$ is not uniquely determined. However we 
simply choose one such lattice for our purpose 
and denote it by $D \ ( = D ( \Sigma , L ) ) $. Then we 
have the exact sequence
of abelian groups:
\be
0 \longrightarrow {\bf n} \stackrel{\iota}{\longrightarrow} 
D \stackrel{\beta}{\longrightarrow} 
L \longrightarrow 0 \ ,
\ele(seqn)
and its dual sequence:
\be
0 \longrightarrow L^* \stackrel{\beta^*}{\longrightarrow} D^*
\stackrel{\imath^*}{\longrightarrow} {\bf n}^* \longrightarrow  0   \ .
\ele(seqndual)
Regard (\req(seqn)) and (\req(seqndual)) as the
corresponding 1-parameter subgroups and characters respectively
for the exact sequence of tori:
\[
0 \longrightarrow \TZ ( {\bf n} ) \longrightarrow \TZ ( D) 
\longrightarrow \TZ ( L ) \longrightarrow 0 \ ,
\]
and identify the following  vector 
spaces:
\[
D_{\KZ} = ( D_0 )_{\KZ} = \bigoplus_{ \delta \in \Sigma^{1} } 
{\KZ} e_{\delta} \ , \ \ \mbox{for \ } \KZ = \QZ , \RZ , \CZ
\  .
\]
Denote 
\[
 \{ e_{\delta}^* \}_{ \delta \in \Sigma^{1}} = \mbox{ \ 
the \ basis \ of \ } D^*_{\KZ}  \mbox{ \ dual \ to \ } 
\{ e_{\delta} \}_{\delta \in \Sigma^{1}} . 
\]
For an element $\sigma$ in a fan $\Sigma$ in $L_{\RZ}$,
let $\tilde{ \sigma}$ be the simplicial cone in $D_{ \RZ }$
defined by
\be
\widetilde{\sigma}  : = \sum_{\delta \in \sigma \cap \Sigma^{1}}
 \RZ_{\geq 0} e_{\delta}  \ \subset \ D_{\RZ} \ ,
\ele(sigmatilde)
and $\widetilde{ \Sigma}$ be the collection of all
such $\widetilde{\sigma}$, 
\[
\widetilde{ \Sigma} : = \{ \widetilde{\sigma} \ | \ \sigma \in
\Sigma \ \} \ .
\] 
Then $\widetilde{ \Sigma}$ is a fan in $D_{ \RZ }$ with its support 
$| \widetilde{ \Sigma } |$ contained in the 
boundary of the first quardrant cone $\Omega$ in $D_{\RZ}$:
\[
\Omega := \sum_{\delta \in \Sigma^{1} } \RZ_{\geq 0} e_{\delta} \ .
\]
Let ${\CZ}_{ ( \Sigma , L ) }$ be the $\TZ ( D)$-variety 
assocoated to the fan $\widetilde{\Sigma }$. It is an open 
subset of its closure $\overline{{\CZ}_{ ( \Sigma , L ) }}$, 
which is a $\TZ ( D)$-affine variety: 
\[
{\CZ}_{ ( \Sigma , L ) } \subset \ 
\overline{{\CZ}_{ ( \Sigma , L ) }}
  = \mbox{Spec } \CZ [ \Omega^* \cap D^* ]  .
\]
The above spaces can be realized as orbifolds in the following manner.
Associate to the pair $( \Sigma , L_0 )$, one 
has the $\TZ ( D_0 )$-toric variety 
$\overline{{\CZ}_{ ( \Sigma , L_0 ) }}$, which is isomorphic 
to $\CZ^d$:
\[
\overline{{\CZ}_{ ( \Sigma , L_0 ) }} = \bigoplus_{\delta \in 
\Sigma^{1} } \CZ e_{\delta} \simeq \CZ^d \ ,
\]
and 
\[
{\CZ}_{ ( \Sigma , L_0 ) } = 
\overline{{\CZ}_{ ( \Sigma , L_0 ) }} - 
 \bigcup_I \{ \sum z_{\delta} e_{\delta}
\ | \ z_{\delta} = 0 \mbox{ for } \ \delta \in I  \} \ ,
\]
where the index $I$ runs over subsets of $\Sigma^{1}$ 
with $I \neq \sigma \cap \Sigma^{1}$ for some 
$\sigma \in \Sigma$.
Since the torus $\TZ ( D_0 )$ is a finite cover of 
$\TZ ( D )$ with the kernel $D/ D_0$:
\[
0 \longrightarrow D/ D_0 \longrightarrow \TZ ( D_0 )
\longrightarrow \TZ ( D ) \longrightarrow 0 \ ,
\]
one has the isomorphisms:
$$
\begin{array}{lll}
\overline{{\CZ}_{ ( \Sigma , L ) }} & \simeq 
\overline{{\CZ}_{ ( \Sigma , L_0 ) }} / ( D/ D_0 ) &\simeq 
\CZ^d / ( D/ D_0 ) \ , \\ [ 2mm ]
{\CZ}_{ ( \Sigma , L ) } & \simeq 
{\CZ}_{ ( \Sigma , L_0 ) } / ( D/ D_0 ) \ .&
\end{array}
$$
Therefore $\overline{{\CZ}_{ ( \Sigma , L ) }} - 
{\CZ}_{ ( \Sigma , L ) }$ 
is an union of affine subvarieties of 
$\overline{{\CZ}_{ ( \Sigma , L ) }}$ having the codimension 
$\geq 2$. Since    
${\CZ}_{ ( \Sigma , L ) }$ is an union of ${\TZ}(D)$-orbits, 
it is stable under the action of ${\TZ}({\bf n})$.  
The linear map $\beta$ of (\req(defbeta)) induces a map from 
$\widetilde{\Sigma }$ to $\Sigma$ sending 
$\widetilde{\sigma}$ to $\sigma $, hence a morphism:
\be
\pi : {\CZ}_{ ( \Sigma , L ) } \longrightarrow 
{\PZ}_{ ( \Sigma , L )} \ .
\ele(principlesigma)
Composing $\pi$ with the finite projection from 
${\CZ}_{ ( \Sigma , L_0 ) }$ to ${\CZ}_{ ( \Sigma , L ) }$, 
one obtains the morphism 
\be
\pi_0 : {\CZ}_{ ( \Sigma , L_0 ) } \longrightarrow 
{\PZ}_{ ( \Sigma , L )} \ .
\ele(sigma0)
The coordinates of the affine spaces 
${\CZ}_{ ( \Sigma , L_0 ) }$ or ${\CZ}_{ ( \Sigma , L ) }$
can be regarded as the generalized homogeneous coordinates for the toric
variety $\PZ_{( \Sigma, L)}$ in \cite{Au} \cite{Cox} ( see also
\cite{Gui} \cite{O94} ). And ${\PZ}_{( \Sigma , L )}$ can be considered
as the set of
${\TZ} ({\bf n})$-orbits in ${\CZ}_{ ( \Sigma , L ) }$:
\[
{\CZ}_{ ( \Sigma , L ) } / {\TZ} ({\bf n}) = {\PZ}_{ ( \Sigma , L )} \ .
\]
In fact, ${\CZ}_{ ( \Sigma , L ) }$ is a principal 
${\TZ} ({\bf n})$-bundle over ${\PZ}_{ ( \Sigma , L )}$ in the sense
of Definition 1.
\par \vspace{0.2in} \noindent
{\bf Lemma 2. } \ With the ${\TZ}({\bf n})$-action on 
${\CZ}_{ ( \Sigma , L ) } $, the 
fibration (\req(principlesigma)) is a principle orbifold
${\TZ}({\bf n})$-bundle. 
\par \vspace{0.2in} \noindent
{\it Proof. } It is known that 
$U_{\sigma} \ ( = \mbox{Spec} \ {\CZ} [ \ \check{ \sigma } 
\cap L^* ] )  , \sigma \in \Sigma$ , form an open affine cover
of  
$ {\PZ}_{ ( \Sigma , L )} $ , and the
${\TZ}(L)$-orbits in ${\PZ}_{\Sigma}$ are in an 
one-to-one correspondence with the open cones 
$\stackrel{\circ}{\sigma}$:
\[
\stackrel{\circ}{\sigma} \longleftrightarrow \mbox{orb} 
( \sigma ) \ , \ \
\sigma \in \Sigma \ .
\]
Over $U_{\sigma}$, $\pi^{-1} ( U_{\sigma} )$ is the affine 
open subvariety $V_{\widetilde{\sigma}} $ of 
$\CZ_{( \Sigma , L )}$,
\[
V_{\widetilde{\sigma}} : = \mbox{Spec} \ {\CZ} 
[ \ \check{ \widetilde{\sigma} } \cap D^* ] , 
\]
and the restriction of $\pi$,
\be
\pi_{\sigma} : V_{\widetilde{\sigma}} \longrightarrow
U_{\sigma} \ ,
\ele(VU)
is the
morphism of toric varieties induced by $\beta$.
For $x \in $ orb ($\sigma) $, 
$U_{\sigma}$ is an open neighborhood of $x$. 
Let $\{ \delta_j \}_{j = 1}^m $ be the intersection of
$\Sigma^{1}$ with $\sigma$, and 
$L_{\sigma}$ denotes the sublattice
$L \cap \sum_{j=1}^m \RZ \delta_j$ of $L$. Let 
$\rho_k , ( 1 \leq k \leq n-m ),$ be elements in $L$ 
which form a basis of $L$ by adding some basis-elements of 
$L_{\sigma}$. Consider the following lattices:
$$
\begin{array}{ll}
L^\prime = & \sum_{j=1}^m \ZZ \delta_j + 
\sum_{k=1}^{n-m} \ZZ \rho_k \ , \\ [2mm]
D^\prime = & \beta^{-1} ( L^\prime ) \ .
\end{array}
$$
Then the following sequence is exact:
\[
0 \longrightarrow {\bf n} \stackrel{\iota}{\longrightarrow} 
D^\prime \stackrel{\beta}{\longrightarrow} 
L^\prime \longrightarrow 0 \ .
\]
Through the homomorphism $\beta$,  $D / D^\prime$ is isomorphic to 
$L / L^\prime$. Denote the group
\be
G_x : = D / D^\prime \simeq L / L^\prime \ .
\ele(Gx)
One has an exact sequence of tori:
\[
0 \longrightarrow \TZ ( {\bf n} ) \longrightarrow \TZ ( D^\prime) 
\longrightarrow \TZ ( L^\prime ) \longrightarrow 0 \ ,
\]
and $G_x$ is imbedded into $\TZ (D^\prime)$ 
and $\TZ (L^\prime )$ compatibly with the above epimorphism.
Denote
\[
U^\prime_{\sigma} = \mbox{Spec} \ \CZ [ \check{ \sigma } 
\cap L^{\prime *} ] \ , \ \
V^\prime_{\widetilde{\sigma}} = \mbox{Spec} \ \CZ [ \check{ \widetilde{\sigma} } 
\cap D^{\prime *} ] \ . 
\]
Then $\beta$ induces an equivariant morphism 
(with respective to $\TZ (D^\prime)$ and $\TZ (L^\prime)$),
\[
\pi^\prime_{\sigma} : V^\prime_{\widetilde{\sigma}} \longrightarrow 
U^\prime_{\sigma} \ .
\]
$\pi^\prime_{\sigma}$ defines a principle 
$\TZ ( {\bf n} )-$bundle on which the group $G_x$ acts 
as a bundle morphism through (\req(Gx)). It is easy to see that
$\pi^\prime_{\sigma}$ is equivalent to the trivial bundle,
\[
pr \ :  U^\prime_{\sigma} \times \TZ ( {\bf n} ) \longrightarrow 
U^\prime_{\sigma} \ \ .
\]
Since the $G_x$-quotient of $\pi^\prime_{\sigma}$ gives the filbration
(\req(VU)), therefore the conditions of Defintion 1 for
the fibration (\req(principlesigma)) are satisfied.
$\Box$ \par \vspace{0.2in} \noindent

\section{ Rational Picard groups of toric varieties }
Let $\PZ_{( \Sigma, L)}$ be the same as in the previous section.
We are going to describe the 
rational Picard group of $\PZ_{( \Sigma , L )}$
in terms of orbifold line bundles. 
Since 
\[
H^q ( \PZ_{( \Sigma , L )} , {\cal O} ) = 0 \ , \ \ \mbox{for} \ q \geq 1 \ ,
\]
we have a natural identification:
\be
\mbox{ Pic} ( \PZ_{( \Sigma , L )})_{ \KZ } \simeq 
H^2 ( \PZ_{( \Sigma , L )} , \KZ )
\ , \ \ \ \KZ = \QZ , \RZ , \CZ \ ,
\ele(PicH2)
by assigning to a line bundle its Chern class.
It is known that the rational cohomology ring of $\PZ_{( \Sigma , L )}$
is the $\QZ$-algebra 
generated by the classes of toric divisors \cite{Dan} 
\cite{F} \cite{O88}. Hence we have the following result,  
(see \cite{O88} Corollary 2.5.): 
\par \vspace{0.2in} \noindent
{\bf Lemma 3. } \ Pic$ ( \PZ_{( \Sigma , L )} )_{\QZ}$ is 
the $\QZ$-space generated 
by the "line bundles" associated to the toric
divisors in $\PZ_{( \Sigma , L )}$.
$\Box$ \par  \noindent   
(Here a non-zero multiple of a toric divisor $D$ is a Cartier divisor, 
so one can talk about its "line bundle" in the rational 
Picard group of $\PZ_{( \Sigma , L )}$    ).

By Lemma 1 and 2, there associates an orbifold line bundle ${\cal O} ( \chi ) $ over
${\PZ}_{( \Sigma , L )}$ from the
$\TZ ( {\bf n} )$-orbifold bundle ${\CZ}_{ ( \Sigma , L ) } $
for a character $\chi$ of 
${\TZ} ({\bf n})$. Hence we have the homomorphism
\be
{\bf n }_{ \QZ }^* \ ( \  = 
{\bf n } ( \Sigma , L )_{ \QZ }^* \ ) \  \longrightarrow \mbox{Pic} 
( {\PZ}_{( \Sigma , L )} )_{ \QZ } \ ,
\ele(nPic)
which sends an element $\chi$ of ${\bf n
}^*$ to ${\cal O} ( \chi ) $.
Define the group
\[
G ( \Sigma , L ) : \ = \{ \varphi : {\CZ}_{ ( \Sigma , L ) }
\longrightarrow {\CZ}_{ ( \Sigma , L ) } ,  \ 
\TZ ( {\bf n} )-\mbox{equivariant \ biregular \ morphism 
\ } \} .
\]
Note that each element in $G ( \Sigma , L )$ can always be extended to
an automorphism of $\overline{{\CZ}_{ ( \Sigma , L ) }}$. 
It has been shown in \cite{Cox} that $G ( \Sigma , L )$ is an 
affine algebraic group, consisting of all the
$\TZ ( {\bf n} )$-biregular morphism of
$\overline{{\CZ}_{ ( \Sigma , L ) }}$. 
The $G ( \Sigma , L )$-action on ${\CZ}_{ ( \Sigma , L ) }$ 
induces a bundle morphism:
$$
\begin{array}{ccl}
G ( \Sigma , L ) \times {\CZ}_{ ( \Sigma , L ) } & \longrightarrow
 & {\CZ}_{ ( \Sigma , L ) } \\
 \downarrow \ \ \ \ \ & & \downarrow  \\  
G ( \Sigma , L ) \times \PZ_{ ( \Sigma , L)} & 
\longrightarrow & \PZ_{ ( \Sigma , L)}  \ ,
\end{array}
$$
hence a $G ( \Sigma , L )$-action on ${\cal O} ( \chi )$. Therefore
we obtain a representation of $G ( \Sigma , L )$ on the vector space
$H^0 ( \PZ_{ ( \Sigma , L)} , {\cal O} ( \chi ) )$.
\par \vspace{0.2in} \noindent
{\bf Theorem 1. } \ For a $n$-dimensional lattice $L$ and a complete 
fan $\Sigma$ in $L_{\RZ}$, let $d$ be the number of irreducible 
toric divisors
in $\PZ_{( \Sigma , L )}$  and
${\bf n}$ be the 
$(d - n)$-dimensional lattice defined by (\req(defn)). Then 
$\mbox{Pic } ( {\PZ}_{( \Sigma , L )} )_{\QZ}$ is a $(d - n)$-dimensional vector 
space generated by toric divisors, and the homomorphism 
(\req(nPic)) induces a natural identification: 
\[
\mbox{Pic }(  {\PZ}_{( \Sigma , L )} )_{\QZ} = 
{\bf n}^*_{\QZ} \ ,
\]
under which a character $\chi$ in ${\bf n}^*$ corresponds 
to the orbifold line bundle ${\cal O} ( \chi )$ over ${\PZ}_{( \Sigma , L )}
$. In particular, the "line bundle" associated to the 
linear equivalent class of 
toric divisor $e_{\delta}$ in $\PZ_{( \Sigma , L )}$ is 
given by $\iota^* (  e_{\delta}^* )$ under the above identification.  
\par \vspace{0.2in} \noindent
{\it Proof. } We are going to show 
the homomorphism (\req(nPic)) is bijective.
The coordinates of ${\CZ}_{ ( \Sigma , L ) }$ are determined 
by the basis $ \{ e_{\delta}^* \}_{ \delta \in \Sigma^{1}} $ 
of $D^*$, and the function
of ${\CZ}_{ ( \Sigma , L ) }$ for $e_{\delta}^*$ is a global
section of ${\cal O} ( \iota^* ( e_{\delta}^* ) )$ with
its zeros at the toric divisor $e_{\delta}$. Hence the 
surjectivity of (\req(nPic)) follows from Lemma 3.
Assume $\chi$ is a character of 
${\TZ} ({\bf n})$ such that
${\cal O} ( \chi^N ) $ is the trivial line bundle over
${\PZ}_{( \Sigma , L )}$ for some positive integer $N$.  
There exists a never-vanishing holomorphic function $f$ on 
${\CZ}_{ ( \Sigma , L ) }$, which extends to $\overline{
{\CZ}_{ ( \Sigma , L ) }}$ by Hartogs' theorem,    
with 
\[
f ( z t ) = \chi^{N} ( t ) f ( z ) \ , \ \ \mbox{for} \ \ t \in {\TZ} ({\bf n})
, \ \ z \in \overline{{\CZ}_{ ( \Sigma , L ) }}\ .
\]
This implies
\[
\chi^{N} ( t ) f ( 0 ) = f ( 0 ) \neq 0  , \ \ \mbox{for} \ \ t \in {\TZ} ({\bf n})
,
\]
here $0$ is the 0-dimensional $\TZ ( D )$-orbit in
$\overline{{\CZ}_{ ( \Sigma , L ) }}$.
Therefore $\chi^N$ is the trivial character, and  
$\chi = 0 \in {\bf n}^*_{\QZ}$.
This shows the injectivity of (\req(nPic)).
$\Box$ \par \vspace{0.2in} \noindent
{\bf Remark. } For the rest of this paper, we shall use the additive 
operation, instead of the
multiplicative one, for the character group of 
$\TZ ( { \bf n } )$ which will be identified with 
${ \bf n }^*$. Then the canonical sheaf of 
$\PZ_{ ( \Sigma , L )}$ is given by the element 
$\iota^* ( \kappa ) $
in ${\bf n}^*$,  where 
\be
\kappa  : = - \sum_{ \delta \in \Sigma^{1}}
e_{\delta}^* \ \in  D_0^* \ \subset D_{\QZ}^* \ ,
\ele(defkappa)
e.g. for the ordinary projective $n$-space $\PZ^n$, the 
$e_{\delta}^*$'s form the $(n+1)$ homogeneous coordinates 
$x_1, \ldots, x_{n+1}$ and  $\iota^* ( \kappa ) $ is equal to 
the line bundle associated to $- \sum_{i=1}^{n+1}D_i $ 
where $D_i$ the divisor of $\PZ^n$ defined by $x_i = 0$.

Now we describe some properties related 
to the ampleness of ${\cal O} ( \iota^* \rho )$ for $\rho \in 
D^*$. For a $n$-dimensional lattice $L$, define 
the $(n+1)$-dimensional lattice,
\[
\overline{ L } =  L \bigoplus \ZZ \ ,
\]
and the element, 
\[
q_L : = ( 0, 1 ) \in \overline{ L } \ .
\]
By the natural inclusion and projection, one 
has the splitting exact sequence:
\be
0 \longrightarrow  \ZZ \longrightarrow 
\overline{L} \longrightarrow  
 L \longrightarrow 0  \ .
\ele(Lseq)
It corresponds to the exact sequence of tori:
\be
0 \longrightarrow  \CZ^* \longrightarrow 
\TZ ( \overline{L}) \longrightarrow  
\TZ ( L ) \longrightarrow 0  \ ,
\ele(seqnLbar)
whose characters are given by the dual sequence of (\req(Lseq)):
\[
0 \longrightarrow  \ L^* \longrightarrow 
\overline{L^*} \longrightarrow  
 \ZZ \longrightarrow 0  \ .
\]
Let $\rho$ be an element in $D_{\QZ}^*$ with the property
\[
\rho_{\delta} : = \rho ( e_{\delta} ) \in \ZZ \ , \ \ \
\mbox{for \ } \delta \in \Sigma^{1} \ .
\]
The piecewise linear functional $f_{\rho}$ on $ L_{\RZ} ( = | \Sigma | )$,
\[
f_{\rho} : L_{\RZ} \longrightarrow \RZ \ ,
\]
is defined by $f_{\rho} ( \delta ) = 
\rho_{\delta}$ for $\delta \in \Sigma^{1}$, and 
linear on each cone in $\Sigma$ . Consider the graph of
$f_{\rho}$:
\[
 L_{\RZ} \longrightarrow \ \overline{L}_{\RZ} 
\ ( =  L_{\RZ}  \bigoplus \RZ ) \ , \ \ \ 
x \mapsto ( x , f_{\rho} ( x )  ) \ .
\]
The fan $\Sigma$ of $L_{\RZ}$ gives rise to a simplicial cone decomposition
of the above graph, denoted by
\be
\Sigma ( \rho ) = \{ \gamma ( \sigma ) \ | \ \sigma \in \Sigma \ \} 
\ , \ \gamma ( \sigma )
: = \mbox{the\ graph\ of\ } f_{\rho} \mbox{ \ over\ the\ region 
\ } \sigma \ .
\ele(Sigmarho)
Then $\Sigma (\rho)$ is a rational  
fan in $\overline{L}_{\RZ}$ with its
support $| \Sigma ( \rho ) |$ equal to the graph of $f_{\rho}$.
Let $E^*_{ \rho }$ be the $\TZ ( \overline{L} )$-toric variety 
associated to $ \Sigma ( \rho ) $, hence it is acted by 
$\CZ^*$ via the sequence (\req(seqnLbar)).
One has the equivariant morphism from $\TZ ( \overline{L} )$-variety
$E^*_{ \rho }$ onto the $\TZ (L)$-variety $\PZ_{ ( \Sigma , L)}$,
\[
\pi_{\rho} : E^*_{ \rho } \longrightarrow \PZ_{ ( \Sigma , L)} \ ,
\]
extending the projection of $\TZ ( \overline{L} )$ to $\TZ
(L)$. It is not hard to see that the above fibration $\pi_{\rho}$ is a 
principal orbifold $\CZ^*$-bundle. We are going to construct 
a partial 
compactification of $E^*_{ \rho }$ along fibers of
$\pi_{\rho}$ as follows. Let ${\cal C}_{ \rho }$ be the
epigraph of $f_{\rho}$, i.e. 
\be
{\cal C}_{ \rho } : = \{ ( x , \alpha ) \in \overline{L}_{\RZ} \ | \
\alpha \geq f_{\rho} ( x ) , \ x \in L_{\RZ} \} \ .
\ele(Crho)
${\cal C}_{ \rho }$ is a $(n+1)$-dimensional cone whose 
interior contains the element $q_L$. 
A simplicial cone decomposition of ${\cal C}_{ \rho }$ is
obtained by starring $\Sigma ( \rho ) $ at $q_L$:
\be
\Sigma ( \rho)*q_L = \bigcup_{\gamma ( \sigma ) \in 
\Sigma ( \rho )} \{ \mbox{ faces \ of } \ ( \gamma ( \sigma )
+ \RZ_{\geq 0} q_L ) \  \} \ .
\ele(Sigtil)
which is a fan in $\overline{L}_{\RZ}$ containing 
$\Sigma ( \rho )$ with
$| \Sigma ( \rho)*q_L | = {\cal C}_{\rho}$. Let 
$E_{\rho}$ be the 
$\TZ ( \overline{L} )$-variety associated to the fan 
$\Sigma ( \rho)*q_L$.
The projection, 
\[
\Sigma ( \rho)*q_L \longrightarrow \Sigma  ,
\]
gives rise to an equivariant morphism from 
$E_{ \rho }$ onto $\PZ_{ ( \Sigma , L)}$.  
By the commutative diagram 
\bea(ccl)
E^*_{ \rho } \ \ & \subset & E_{ \rho } \\
\downarrow \pi_{\rho}  & & \downarrow \pi_{\rho} \\ [2mm] 
\PZ_{ ( \Sigma , L)} & = & \PZ_{ ( \Sigma , L)}  \ ,
\elea(E*E)
$E_{ \rho } - E^*_{ \rho }$ isomorphic to $\PZ_{ ( \Sigma , L)}$
under the map $\pi_{\rho}$, and $E_{\rho}$ is a partial 
compactification of $E^*_{ \rho }$. On the other hand,
the homomorphism,
\be
p_{\rho} : D \longrightarrow \overline{L} \ \ , \ \ 
 x \mapsto (  \beta ( x ) , \rho ( x )  \ ) \ ,
\ele(rho)
induces a tori-morphism
\[
\TZ ( D ) \longrightarrow \TZ ( \overline{L} ) .
\] 
One has the commutative diagram of exact sequences:
\bea(llll)
0 \longrightarrow & {\bf n} \stackrel{\iota}{\longrightarrow} 
& D \stackrel{\beta}{\longrightarrow} & L 
\longrightarrow 0 \\ [2mm]
& \downarrow \iota^*{\rho} & \downarrow p_{\rho}
& \|   \\ [1mm]
0 \longrightarrow & \ZZ \longrightarrow 
&\overline{L} \longrightarrow  
& L \longrightarrow 0  \ ,
\elea(comexseq)
and its dual diagram:
\bea(llll)
0 \longrightarrow & L^* \longrightarrow 
&\overline{L^*} \longrightarrow  
& \ZZ \longrightarrow 0  \\ [2mm]
& \| & \downarrow p^*_{\rho}
& \downarrow \jmath  \\ 
0 \longrightarrow & L^* \stackrel{\beta^*}{\longrightarrow} 
&D^* \stackrel{\iota^*}{\longrightarrow} & {\bf n}^* 
\longrightarrow 0 ,
\elea(comexseqdual)
where $\jmath$ is the linear map sending 1 to $\iota^*{\rho}$, 
and $p^*_{\rho}$ has the expression 
\[
p^*_{\rho} : \overline{L^*} \longrightarrow D^* \ , \ \
( l^\prime , a  ) \mapsto  \beta^* ( l^\prime ) + a \rho \ \ \
 \mbox{for} \ a \in \RZ , \ l^\prime \in L^*_{\RZ} \ .
\]
For $\sigma \in \Sigma$, $\gamma ( \sigma ) $ is generated 
by $\{ p_{\rho} ( e_{\delta} ) \ | 
\delta \in \sigma \cap \Sigma^{1} \}$, and the linear map,
\[
p_{\rho} : D_{\RZ} \longrightarrow \overline{L}_{\RZ} \ ,
\]
sends the simplicial cone $\widetilde{\sigma}$ of
(\req(sigmatilde)) in $D_{\RZ}$ to $\gamma ( \sigma ) $.
Hence we have the expression  
\[
\Sigma ( \rho )  = \{ p_{\rho} ( \widetilde{\sigma} ) \ | \
\sigma \in \Sigma  \} \ ,
\]
from which one obtains the induced "bundle" morphism: 
$$
\begin{array}{ccc}
{\CZ}_{ ( \Sigma , L ) }  & \longrightarrow  & E^*_{ \rho } \\ [2mm]
 \downarrow \pi &  
& \downarrow \pi_{\rho}  \\ [2mm]
\PZ_{ ( \Sigma , L)} & = & \PZ_{ ( \Sigma , L)} .
\end{array}
$$
With the character $\iota^* ( \rho )$ of $\TZ ( {\bf n })$ on 
${\CZ}_{ ( \Sigma , L ) }$ and the identity character of
$\CZ^*$ on $E^*_{ \rho }$ in Definition 2, one has the commutative diagram: 
$$
\begin{array}{ccc}
{\cal O} ( \iota^*(\rho) ) = {\CZ}_{ ( \Sigma , L ) } \times_{\TZ ( {\bf n } )} \CZ  
& \simeq  & E^*_{ \rho } \times_{\CZ^*} \CZ = E_{ \rho }  \\ [2mm]
 \downarrow  &  
& \downarrow   \\ [2mm]
\PZ_{ ( \Sigma , L)} & = & \PZ_{ ( \Sigma , L)} .
\end{array}
$$ 
In this manner, $E^*_{ \rho }$ is identified with the
complement of the zero section of ${\cal O} ( 
- \iota^* \rho )$ over $\PZ_{ ( \Sigma , L)}$, and  
(\req(E*E)) is equivalent to the diagram:
$$
\begin{array}{ccl}
{\cal O} ( - \iota^* \rho )- \mbox{(zero-section)}  & \subset & {\cal O} ( - \iota^* \rho ) \\
 \downarrow  & & \downarrow \\ 
\PZ_{ ( \Sigma , L)} & = & \PZ_{ ( \Sigma , L)}  \ .
\end{array}
$$

From now on, we consider only those ${\cal O} ( \iota^* \rho )$
having non-trivial sections. Without loss of 
generality we shall always make the following assumption
on $\rho$ for the rest of this paper: 
\be
\rho = \sum_{ \delta \in \Sigma^1 } \rho_{\delta} \ e^*_{\delta}  
\in D^* - \{  0  \} \ , \ \ \
\rho_{\delta} \in \ZZ_{\geq 0} \ \ \mbox{for \ all \ } \delta \ ,
\ele(rhopost)
or equivalently, $f_{\rho}$ is a non-negative and non-trivial function on 
$L_{\RZ}$. Note that the cone ${\cal C}_{ \rho }$ 
in $\overline{L}_{\RZ}$ is non-convex in 
general; while the $p_{\rho}$-image of first quardrant 
cone $\Omega$ of $D_{\RZ}$, $p_{\rho} ( \Omega )$, is always 
a convex polyhedral cone, which is characterized as the convex cone in $\overline{L}_{\RZ}$
generated by $( \delta , \rho_{\delta} ) , \delta \in \Sigma^{1}$. Since $\rho$ is not trivial,
$p_{\rho} ( \Omega )$ is a $(n+1)$-dimensional cone 
containing $q_L$ and the relation 
\[
{\cal C}_{ \rho } \subseteq p_{\rho} ( \Omega ) 
\]
holds. 
The equality of the above sets is related to the convexity 
of $\rho$ in the context of ample toric line bundles  
\cite{Dan} \cite{KKMS} \cite{O78}.
As the generating set $\sigma \cap \Sigma^{1}$ of a 
$n$-dimensional cone $\sigma$ in $\Sigma$ forms a 
$\QZ$-basis of $L_{\QZ}$, 
there is an unique $\QZ$-linear functional $l_{\sigma}$ of $L_{\QZ}$ 
with the property :
\[
l_{\sigma} ( \delta ) = \rho_{\delta} \ , \ \mbox{for} \ \ \
\delta \in \sigma \cap \Sigma^{1} \ .
\]
The following definition is given in \cite{O88} \cite{OP}: 
\par \vspace{0.2in} \noindent
{\bf Definition 3. } For an element $\rho$ in 
(\req(rhopost)), we say

(i) $\rho$ is convex if the 
functional $f_{\rho}$ on $L_{\RZ} $ is a convex function,
i.e. 
$$
\begin{array}{lcl}
& l_{\sigma} ( x ) \leq f_{\rho} ( x ) \ , & 
 \mbox{ \ 
for \ } \sigma \in \Sigma^{(n)} , \ x \in L_{\RZ} , \\
\Longleftrightarrow &  
f_{\rho} ( x + y ) \leq f_{\rho} ( x ) + f_{\rho} ( y ) , &
\mbox{ \ for \ } x , y \in L_{\RZ} \ . 
\end{array}
$$ 
 
(ii) $\rho$ is strictly convex if $f_{\rho}$ is a strictly convex 
function, i.e. for $\sigma \in \Sigma^{(n)}$ and 
$x \in L_{\RZ}$,
$$ 
\begin{array}{lll}
l_{\sigma} ( x ) \leq f_{\rho} ( x ) \ , & 
\mbox{and} & \mbox{"=" holds \ iff } \ x \in \sigma .
\end{array}
$$
\par \vspace{0.2in} \noindent
{\bf Proposition 1. } \ The following conditions are equivalent:

(i) $\rho$ is convex if and only if ${\cal C}_{\rho}$ is a convex 
cone. In this situation, the following equality holds:
\[
{\cal C}_{\rho} =  p_{\rho} ( \Omega ) .
\]

(ii) $\rho$ is strictly convex if and only if ${\cal C}_{\rho}$
is a convex cone in $\overline{L}_{\RZ}$ having all the 
proper faces to be simplicial and with  
$\{ (  \delta , \rho_\delta  ) \ | \ 
\delta \in \Sigma^{1} \}$ as the fundamental generators
, (i.e. irredundant ones). In this situation, the cone 
${\cal C}_{\rho}$
is strongly convex, i.e.   
${\cal C}_{\rho} \cap ( - {\cal C}_{\rho} ) = \{ 0 \}$ .
\par \vspace{0.2in} \noindent
{\it Proof. } Since the elements 
$(  \delta , \rho_{\delta} )$ of ${\cal C}_{\rho}$ 
generate the convex cone $p_{\rho} ( \Omega )$,
the convex condition for ${\cal C}_{\rho}$ is equivalent to the 
equality of ${\cal C}_{\rho}$ and $p_{\rho} ( \Omega )$. 
For $\sigma \in \Sigma^{(n)}$, $\gamma ( \sigma )$ is a  
$n$-dimensional cone in $\overline{L}_{\RZ}$
generated by elements $
(  \delta , \rho_{\delta}  ) $ 
for $ \delta \in \sigma \cap \Sigma^{1}$. 
The graph of $l_{\sigma}$ in $\overline{L}_{\RZ}$ coincides
with $\gamma ( \sigma )$ over the region $\sigma$. Therefore the 
condition of $\gamma ( \sigma ) \subseteq \partial p_{\rho} ( \Omega )$ becomes
$\rho ( x ) \geq l_{\sigma} ( x )$ for $x \in  \Omega $, which   
is equivalent to
\[
f_{\rho} ( \delta ) \ ( = \rho_{\delta}  ) \ \geq 
l_{\sigma} (\delta ) \ \ 
\mbox{for} \ \delta \in \Sigma^{1} \ . 
\]
Therefore we obtain (i). It is not hard to see that 
the stricly convexity of $\rho$ is equivalent to
the irredundant condition of the generators $(  \delta ,
 \rho_{\delta} )$'s together with the simplical proper face-property 
of $p_{\rho} ( \Omega )$. 
In this situation, we have ${\cal C}_{\rho} =  p_{\rho} ( \Omega ) $. 
If $p_{\rho} ( \Omega )$ contains a linear space generated by a 
non-trivial element $v$ in $\overline{L}_{\RZ}$, write
\[
v = (  x_v , \alpha ) \ , \ x_v \in L_{\RZ} \ , \
\ \alpha \in \RZ \ .
\]
$x_v$ belongs to some $n$-dimensional cone $\sigma $ of $\Sigma$. 
Since $p_{\rho} ( \Omega )$ is contained in the "upper" 
half-space of $\overline{L}_{\RZ}$, $\alpha$ is equal to 0 , hence 
$ \pm v \in \partial p_{\rho} ( \Omega )$.
This implies 
$ l_{\sigma} ( x_v ) = 0$. By the strictly convex condition of 
$\rho$ and $-x_v \not\in \sigma$, 
 $f_{\rho} ( - x_v ) > l_{\sigma} (- x_v ) = 0 $, which 
contradicts $-v \ 
 ( = (  -x_v , 0 ) ) \in p_{\rho} ( \Omega )$. 
Therefore $C_{\rho}$ is strongly convex.
$\Box$ \par \vspace{0.2in} \noindent
As a consequence of Proposition 1, the ampleness for a  
multiple of anti-canonical divisor of $\PZ_{( \Sigma , L)}$ 
is now determined by the polytope
$\triangle ( L )$ of (\req(deftria)):
\par \vspace{0.2in} \noindent
{\bf Corollary. } Let $\kappa$ be the element (\req(defkappa))
and $r$ a negative rational number such that $r \kappa \in D^*$. Then
$$ 
\begin{array}{lll}
r \kappa : \mbox{ convex } & \Longleftrightarrow &
\triangle ( L ) : \mbox{ convex } \ ; \\
r \kappa : \mbox{ strictly \ convex } & \Longleftrightarrow &
\triangle ( L ) : \mbox{ convex \ cone \ with \ simplical \
proper \ faces \ and  } \\
& & \ \ \ \ \ \ \ \ \ \mbox{ having \ }
\{ \delta \}_{\delta \in \Sigma^{1}}  
\mbox{ as \ a \ set \ of \ minimal \ generators \ } .
\end{array}
$$ 
$\Box$ \par \vspace{0.2in} \noindent
For a convex element $\rho$, ${\cal C}_{\rho}$ is a 
$(n+1)$-dimensional convex cone in the upper half-space 
of $\overline{L}_{\RZ}$. Let 
$\mbox{v}_{\rho}$
be the maximal linear subspace of $L_{\RZ}$ contained in 
${\cal C}_{\rho}$, and $\mbox{v}_{\rho}^{ \perp }$ 
be its annahilating subspace in $L_{\RZ}^*$. Then 
the dual cone ${\cal C}^*_{\rho}$ is a convex cone in 
$ \mbox{v}_{\rho}^{ \perp } + \RZ q_{L^*}$ with $q_{L^*} \in 
\mbox{Int} ( {\cal C}_{\rho}^* )$:
\[
q_{L^*} \in {\cal C}^*_{\rho} \ \subseteq \ 
\mbox{v}_{\rho}^{ \perp } + \RZ q_{L^*} \ \subseteq  
\ \overline{L^*}_{\RZ} \ .
\]
The relation between ${\cal C}^*_{\rho}$ and the dual 
first quardrant cone $\Omega^* ( = \sum_{\delta 
\in \Sigma^{1} } \RZ_{\geq 0} e^*_{\delta} )$ of $ D^*_{\RZ} $
is given by $p_{\rho}^*$:
\[
p_{\rho}^* \ : \ {\cal C}^*_{\rho}  \ \simeq   \ 
\Omega^* \bigcap p_{\rho}^*( \overline{L^*}_{\RZ} ) \ .
\]
Since the boundary of ${\cal C}_{\rho}$ is 
the graph of $f_{\rho}$, every 1-face of ${\cal C}_{\rho}^*$ must 
intersect the $( q_L=1 )$-hyperplane at one point.
${\cal C}_{\rho}^*$ is determined by its
cross section with the 
$( q_L=1 )$-hyperplane.
Let $\triangle ( L^* )_{\rho}$ be the projection of
this cross section to $L_{\RZ}$:
\be
\triangle ( L^* )_{\rho} : = \{ x \in L^*_{\RZ} \ | \ 
(  x , 1 ) \in {\cal C}_{\rho}^* \}  \ .
\ele(tria*) 
Then $\triangle ( L^* )_{\rho}$ is a compact convex polytope in 
$\mbox{v}_{\rho}^{ \perp }$, hence in $L^*_{\RZ}$. In the next section we shall 
show that 
integral elements in $\triangle ( L^* )_{\rho}$ determine
the sections of ${\cal O} ( \iota^*\rho )$ over 
$\PZ_{ ( \Sigma, L ) }$.

\section{ Quasi-smooth hypersurfaces in toric varieties }
First let us recall the notion of quasi-smoothness of a 
hypersurface in a toric variety introduced by
Batyrev and Cox in \cite{BC}. A hypersurface $X$ in 
$\PZ_{( \Sigma , L)}$ is called to be quasi-smooth if $\pi_0^*(X)$ is 
non-singular in ${\CZ}_{ ( \Sigma , L_0 ) }$ where 
$\pi_0$ is the morphism in (\req(sigma0)). This 
concept is equivalent to "simplicially toroidal pairs" by 
Danilov \cite{D91}. In this section we are going to discuss 
some basic properties of quasi-smooth hypersurfaces defined 
by sections of an orbifold line bundle over 
$\PZ_{( \Sigma , L)}$.
Let $\rho$ be a convex element of $D^*$ with
the cone ${\cal C}_{\rho}$   
in $\overline{L}_{\RZ}$,  its dual convex 
cone ${\cal C}_{\rho}^*$ in $\overline{L^*}_{\RZ}$, and
the convex polytope $\triangle ( L^* )_{\rho}$ in $L^*_{\RZ}$.
\par \vspace{0.2in} \noindent
{\bf Proposition 2.} The vector space $
\Gamma ( \PZ_{( \Sigma , L)} , {\cal O} ( \iota^* \rho ) )$ 
has a basis consisting of elements $
z \in \overline{L^*} \cap {\cal C}_{\rho}^* $ with $
< q_L , z > = 1$, which are in one-to-one correspondence with
elements in $L^* \cap \triangle ( L^* )_{\rho}$, i.e.  
$$
\begin{array}{lll}
\Gamma ( \PZ_{( \Sigma , L )} , {\cal O} ( \iota^* \rho ) ) 
& \simeq & \bigoplus \{ \CZ z \ | \ 
z \in \overline{L^*} \cap {\cal C}_{\rho}^* \ , \
< q_L , z > = 1 \ \} \ , \\
& \simeq & \bigoplus \{ \CZ x \ | \ 
x \in L^* \cap \triangle ( L^* )_{\rho} \ \} \ .
\end{array}
$$
\par \vspace{0.2in} \noindent
{\it Proof.} Consider the fibration $E_{\rho}$ over 
$\PZ_{( \Sigma , L)}$ in (\req(E*E)). It is 
isomorphic to the orbifold line bundle 
${\cal O} ( - \iota^* \rho )$ over $\PZ_{ ( \Sigma , L)}$ 
with the scalar multiplication corresponding to the action of 
one-parameter subgroup $q_L$ on $E_{\rho}$. Hence
the sections in
$\Gamma ( \PZ_{(\Sigma, L)} , {\cal O} ( \iota^* \rho ) )$ 
can be regarded as regular functions of $E_{\rho}$,
linear with respect to the $\CZ^*$-action induced by $q_L$. By 
$| \widetilde{\Sigma} ( \rho ) | = {\cal C}_{\rho}$, 
$\overline{L^*} \cap {\cal C}_{\rho}^*$ determines the
regular functions of $E_{\rho}$. Hence a basis of 
$\Gamma ( \PZ_{(\Sigma, L)} , {\cal O} ( \iota^* \rho ) ) $ 
consists of the elements in 
$\overline{L^*} \cap {\cal C}_{\rho}^*$ whose 
$q_L$-values equal to 1. It is easy to see that they are in one-to-one 
correspondence with elements in
$L^* \cap \triangle ( L^* )_{\rho}$. 
$\Box$ \par \vspace{0.2in} \noindent
By the Proposition 2, a section $s$ in
$\Gamma ( \PZ_{(\Sigma, L)} , {\cal O} ( \iota^* \rho ) ) $ 
has the following expression:
\be
s  =  \sum_{ x \in L^* \cap \triangle ( L^* )_{\rho}} 
\alpha_x x \ ,  \ \ \  \ \mbox{or \ equivalently } \ \ \ \ \
s = \sum_{z \in \overline{L^*} \cap {\cal C}^*_{\rho} , 
<q_L, z> = 1 } \alpha_z z \ .
\ele(fexpr)
One can consider $s$ as a function on the 
$\TZ ( \overline{L} )$-variety $E_{\rho}$ ( 
$ \simeq {\cal O} ( - \iota^* \rho )$ ). 
The Newton polygon of $s$ is defined to be the convex hull in 
$L^*_{\RZ}$ spanned by $x$'s with $\alpha_x \neq 0$, hence is 
contained in $\triangle ( L^* )_{\rho}$.
Let $X$ be the hypersurface of $\PZ_{(\Sigma , L)}$ defined by 
the zeros of $s$.
We are going to describe the local defining equation of 
$X$. For $\sigma \in \Sigma^{(n)}$, 
we have the affine open subset $U_{\sigma}$ of 
$\PZ_{( \Sigma , L )}$:
\[
U_{\sigma} = \mbox{ Spec } \CZ [ \breve{ \sigma } \cap L^* ] 
\subset \PZ_{( \Sigma , L )} \ .
\]
Denote 
\be
\{ \delta_i \}_{ i = 1}^n  = \Sigma^{1} \cap \sigma \ ,
\ele(defdeli) 
and let $I_{\sigma}$ be the finite subgroup of $\CZ^{n *}$ 
defined by the image of $L$ under the following homomorphism:
\[
\sum_{i=1}^n \RZ \delta_i \longrightarrow \CZ^{n *} \ , \ \ \ \  
\sum_{i=1}^n x_i \delta_i \mapsto ( e^{2\pi i x_1} ,
\ldots , e^{2\pi i x_n} ) \ .
\]
Then $I_{\sigma}$  
is isomorphic to $L / \sum_{i=1}^n \ZZ \delta_i $, and 
$U_{\sigma}$ is realized as an orbifold through the 
projection:
\be
\psi : \CZ^n \longrightarrow U_{\sigma} 
= \CZ^n / I_{\sigma} \ ,
\ele(psi)
where the dual basis of $\{ \delta_i \}_{i=1}^n$ gives rise to the 
coordinates $( t_1, \ldots , t_n )$ of $\CZ^n$. By 
Proposition 3. 5 of \cite{BC}, quasi-smooth hypersurfaces 
have the following local descrption.   
$X \bigcap U_{\sigma}$ is quasi-smooth if
 $\psi^* ( X \bigcap U_{\sigma} )$ is a smooth hypersurface in
$\CZ^n$, and we call    
$X$ to be quasi-smooth if $X \bigcap U_{\sigma}$ is quasi-smooth 
for all $\sigma \in \Sigma^{(n)}$. 
Let
\[
\left\{ \begin{array}{ll}
p_i = (  \delta_i , \rho_{\delta_i} ) \ , & \mbox{for} \
1 \leq i \leq n , \\
p_{n+1} = q_L \ &.
\end{array}
\right.
\]
All the $p_j$'s are primitive
elements in $\overline{L}$ and they generate a 
$(n+1)$-simplicial cone in $\widetilde{\Sigma}( \rho )$.
The first $n$ elements $p_1, \ldots , p_n ,$ lie in an unique
$n$-dimensional face of ${\cal C}_\rho$, annihilated by 
a vertex $z ( \sigma )$ of the polytope 
${\cal C}_\rho^* \cap ( q_L = 1 )$. Note that $z ( \sigma )$ may not
be an element in $ \overline{L^*}$ unless 
$\triangle ( L^* )_{\rho} $ is the Newton polygon of $s$. 
The dual elements 
$p_j^* , 1 \leq j \leq (n+1) ,$ are in 
$\overline{L^*}_{\QZ}$ with
\[
< p_i , p_j^* > = \delta_{ij} \ , \ 1 \leq i , j \leq (n+1) .
\]
Hence
\[
p_1^*, \ldots , p_n^* \in L_{\QZ}^* \ , \ \ 
\mbox{and \ } p_{n+1}^* = z ( \sigma ) \ ,
\]
and $\{ p_i^* \}_{i=1}^n$ is the dual basis of 
$\{ \delta_i \}_{i=1}^n$.     
By
\[
{\cal C}_\rho^* \subset \ \sum_{j=1}^{n+1} \RZ_{\geq 0} p_j^* \ , 
\]
for $z \in \overline{L^*} \cap {\cal C}^*_{\rho} $ with
$<q_L, z> = 1 $, one has
\be
y_{\sigma} ( z ) : = z - z ( \sigma ) = \sum_{i=1}^n k ( z )_i p_i^* \ , \ \mbox{with} \
k ( z )_i = < p_i , z > \in \ZZ_{\geq 0} \ .
\ele(kz)
Consider $y_{\sigma} ( z )$ as a regular function of $U_{\sigma}$. 
Through (\req(psi)), one has
\[
\psi^* ( y_{\sigma} ( z ) ) = \prod_{i=1}^n t_i^{k ( z )_i} \ ,
\]
and we obtain the defining equation for
$X \bigcap U_{\sigma} $,
\[
 X \bigcap U_{\sigma} \ \ \ : \ \ \ 
\sum_{z \in \overline{L^*} \cap {\cal C}^*_{\rho} , 
<q_L, z> = 1 }
\alpha_z y_{\sigma} ( z ) = 0 \ .
\]
The hypersurface $\psi^* ( X \bigcap U_{\sigma} )$ in
$\CZ^n$ is now defined by
\be
\eta_{\sigma} \ : \ = \sum_{z \in \overline{L^*} \cap {\cal C}^*_{\rho} , 
<q_L, z> = 1 }
\alpha_z \prod_{i=1}^n t_i^{k ( z )_i} = 0 .
\ele(eqnX) 
Now assume ${\cal C}_\rho$ is   
a strongly convex polyhedral cone in $L_{\RZ}$. One has the 
order reversing bijective map:
$$
\begin{array}{ccc}
\{ \mbox{ face \ of } \ {\cal C}_{\rho} \} &
\longleftrightarrow & \{ \mbox{ face \ of \ }
 {\cal C}_{\rho}^*  \} \\ [2mm]
\tau & \leftrightarrow
& \tau^\prime : = {\cal C}_{\rho}^* \cap \tau^{\bot}
\end{array}
$$
with $ \mbox{dim} ( \tau ) + \mbox{dim} ( \tau^\prime ) = n $.
$\tau^\prime$ is called the dual of $\tau$ in 
${\cal C}_\rho^*$. The dual of a facet, (i.e. codimensional
one face), of ${\cal C}_{\rho}$ is 
the 1-face of ${\cal C}^*_{\rho}$ generated by 
$( v , 1 )$ for some vertex $v$ of $\triangle ( L^* )_{\rho}$,
and we have 
\[
v \in L^* \ \Longleftrightarrow \ ( v , 1 ) \in \overline{L^*} \ .
\]
\par \vspace{0.2in} \noindent
{\bf Lemma 4. } Let $X$ be a hypersurface of 
$\PZ_{(\Sigma, L)}$ defined by zeros of a section 
$s \in \Gamma ( \PZ_{(\Sigma, L)} , {\cal O} ( \iota^* \rho ) ) $.
For a vertex $v$ of 
$\triangle ( L^* )_{\rho}$, let 
$\tau$ be the facet of ${\cal C}_{\rho}$ 
dual to $( v , 1 ) \in {\cal C}_{\rho}^*$. 
If the Newton polygon of $s$ contains $v$, 
then 
\[
X \bigcap \bigcup_{\gamma ( \sigma ) \subset \tau} 
\mbox{orb} (  \sigma  ) = \emptyset \ ,
\]
where $\mbox{orb} (  \sigma  )$ is the 
$\TZ ( L )$-orbit in $\PZ_{(\Sigma, L)}$ and  
$\gamma ( \sigma )$ is the cone (\req(Sigmarho)) in the fan 
$\Sigma ( \rho )$ associated to $\sigma$.
\par \vspace{0.2in} \noindent
{\it Proof.} Regard $s$ as a function of 
$E_{\rho}^*$ ( = ${\cal O}( - \iota^* \rho )$ - (zero-section) ) . 
For $t \in \mbox{Int} ( \tau )$ and $x \in 
L^* \cap \triangle ( L^* )_{\rho}$ with $\alpha_x \neq 0 $ in 
(\req(fexpr)), we have 
\[
< t , ( x , 1 )  > \ \geq 0 \ , \ \ \mbox{and \ "=" \ holds \ iff \ }
x = v \ . 
\]
For $ \gamma ( \sigma ) \subset \tau$ and
an element $e$ of $ E_{\rho}^*$ lying over  
$\mbox{orb} (  \sigma  ) $ , 
we have 
\[
s ( e ) = \alpha_v v ( e ) \neq 0 \ .
\]
Therefore the conclusion follows immediately.
$\Box$ \par \vspace{0.2in} \noindent
{\bf Definition 4. } Let $M$ be a lattice. A convex polytope in 
$M_{\RZ}$ is called integral ( with respect to $M$ ) if all its 
vertices are in $M$. 
$\Box$ \par \vspace{0.2in} \noindent
The following result is obvious.
\par \vspace{0.2in} \noindent
{\bf Lemma 5.} The following conditions are equivalent:

(i) Newton polygon of $s$ is equal to  $\triangle ( L^* )_{\rho}$.

(ii) $\triangle ( L^* )_{\rho}$ is an integral polytope in 
$L^*_{\RZ}$, and the coefficient $\alpha_x$ in (\req(fexpr)) 
is non-zero for each vertex $x$ of 
$\triangle ( L^* )_{\rho}$.
$\Box$ \par \vspace{0.2in} \noindent
{\bf Remark.} For an integral polytope 
$\triangle ( L^* )_{\rho}$, a generic section $s$ 
of $\Gamma ( \PZ_{(\Sigma, L)} , 
{\cal O} ( \iota^* \rho ) ) $ always satisfies the above 
lemma, and it defines a quasi-smooth hypersurface of  
$\PZ_{(\Sigma, L)}$.
\par \vspace{0.2in} \noindent   
For a strictly convex element $\rho$, ${\cal C}_{\rho}$
is a strongly convex polyhedral cone in $\overline{L}_{\RZ}$
by Proposition 1, hence
$\triangle ( L^* )_{\rho}$ is a $n$-dimensional convex polytope in
$L^*_{\RZ}$.
Let $A_{\rho}$ be the affine 
$\TZ ( \overline{L} )$-variety associated to the cone 
${\cal C}_{\rho}$ :
\be
A_{\rho} = \mbox{Spec} \ \CZ 
[ \ \overline{L}^* \cap {\cal C}_{\rho}^* ] \ .
\ele(ASpec)
Since $ \{ ( \delta , \rho_\delta  ) \ | \ 
\delta \in \Sigma^{1} \}$ forms a system of fundamental generators of 
${\cal C}_{\rho}$ , $A_{\rho}$ 
is the $\TZ ( \overline{L} )$-variety obtained by adding one point 
( = the 0-dimensional orbit ) to $E_{\rho}^*$. 
In fact, $A_{\rho}$ is the affine variety obtained by blowing down
the zero section of the orbifold line bundle ${\cal O} 
( - \iota^*\rho ) $ over $\PZ_{(\Sigma, L)}$:
\[
 E_{\rho} = {\cal O} ( - \iota^* \rho ) \longrightarrow 
A_{\rho} \ ,
\]
and the one-parameter subgroup $q_L$
gives the $\CZ^*$-action on $A_{\rho}$:
\[
\CZ^* \times A_{\rho} \longrightarrow A_{\rho} \ , \
( \lambda , z ) \mapsto \lambda \cdot z \ .
\]
By Proposition 2, one has the identification:
\[
\Gamma ( \PZ_{(\Sigma, L)} , {\cal O} ( \iota^* \rho ) )  = 
\{ f : A_{\rho} \longrightarrow \CZ \ | \ f ( \lambda \cdot z ) = 
\lambda f ( z ) \ \mbox{for \ } \lambda \in \CZ^* \ , \ z \in 
A_{\rho} \  \} .
\]
The affine coordinates of $A_{\rho}$ can be considered as the 
"homogeneous" coordinates of $\PZ_{(\Sigma, L)}$ for sections of
${\cal O} ( \iota^* \rho ) $. The Lefschetz-type theorem on 
the cohomology of $\PZ_{(\Sigma, L)}$ and its hypersurface 
was obtained by Danilov and Khovanskii \cite{DK}. The 
following form can be derived by a similar 
argument as in Theorem 2 of \cite{AF}.
We repeat it here for the sake of completeness.
\par \vspace{0.2in} \noindent
{\bf Proposition 3. } Let $\rho$ be a strictly convex element in 
$D^*$, and $X$ be a quasi-smooth hypersurface of 
$\PZ_{(\Sigma, L)}$ defined by a 
section of ${\cal O} ( \iota^* \rho )$. Then the homomorphism
\[
H^i ( \PZ_{(\Sigma, L)} , \CZ ) \longrightarrow H^i ( X , \CZ )
\]
induced by inclusion is an isomorphism for $ i < n-1 $, and injective
for $i = n-1$. 
\par \vspace{0.2in} \noindent
{\it Proof. } Let $U$ be the complement of $X$ in 
$\PZ_{(\Sigma, L)}$. The quasi-smooth property for $X$ and $U$ implies
that $U$ is a rational homology manifold where
Poincar$\acute{e}$ duality holds for the cohomology
(with complex coefficients throughout):
\be
H^i ( \PZ_{(\Sigma, L)} , X  ) \ \simeq \ H_{2n-i} ( U ) \ .
\ele(Poincare)
Since ${\cal O} ( \iota^* \rho )$ is an ample line
bundle over $\PZ_{(\Sigma, L)}$, $U$ is an affine variety, therefore
a Stein space. By the acyclic resolution of the constant sheaf
$\underline{\CZ}_U$:
\[
0 \longrightarrow \underline{\CZ}_U \longrightarrow \Omega^{\cdot}_U 
\ ( = \{ \Omega^0_U \stackrel{d}{\longrightarrow} \Omega^1_U 
\stackrel{d}{\longrightarrow} \cdots \stackrel{d}{\longrightarrow} 
\Omega^n_U \longrightarrow 0 \ \} ) \ ,
\]
one has
\be
H^j ( U  ) = 0 \ \ \ \mbox{for} \ j > n \ ,
\ele(HU=0)
hence
\[
H^i ( \PZ_{(\Sigma, L)} , X  ) \ \simeq \ H_{2n-i} ( U  ) = 0 
\ , \ \ \ \mbox{for} \ j \leq n .
\]
By the exact cohomology sequence:
\be
\cdots \longrightarrow H^i ( \PZ_{(\Sigma, L)} , X ) 
\longrightarrow H^i ( \PZ_{(\Sigma, L)} ) \longrightarrow 
H^i ( X ) \longrightarrow H^{i+1} ( \PZ_{(\Sigma, L)} , X ) 
\longrightarrow \cdots 
\ele(MVseq)
the conclusion follows immediately.
$\Box$ \par \vspace{0.2in} \noindent
{\bf Remark.} A similar statement can be found in 
Theorem 3.7 of \cite{DK}. 
The relation between Hodge groups of $X$ and 
$\PZ_{(\Sigma , L )}$  
can be derived from the above proposition, 
for the argument see Appendix (I) .

\section{Examples}
In this section, some examples are given for the illustration 
of the results obtained in previous sections.

(I) The case for $n=1$. 
We have
$$
\begin{array}{cll}
L & = & \ZZ \ , \\
\Sigma & = & \{ \RZ_{\geq 0} \ , \ \RZ_{\leq 0} \ , \ 0 \ \} \ , \\
\Sigma^{1} & = & \{ \delta_{+} , \delta_{-} \}  \ , \
\delta_{\pm} = \pm 1 \ , \\
\PZ_{(\Sigma , L)} & = & \PZ^1 \ ,
\end{array}
$$
hence 
\[
D =  \ZZ^2 \ \ ( \ = \ZZ e_{\delta_+} + \ZZ e_{\delta_-} ) \ .
\]
The map $\beta$ is given by
\[
\beta : D \longrightarrow L \ , \ 
( m_1 , m_2 ) \mapsto m_1 - m_2 \ ,
\]
with the kernel
\[
{\bf n} = \ZZ ( 1 , 1 ) \ .
\]
The $\CZ^*$-bundle (\req(principlesigma)) is the well-known Hopf
fibration:
\[
\pi : \CZ^2 -\{ 0 \} \longrightarrow \PZ^1 \ .
\]
For a given
\[
\rho = ( k, l  ) \in D^* - \{ 0 \} \ , \ 
\  \  k , l  \in \ZZ_{\geq 0} \ ,
\]
one has
\[
{\cal C}_\rho  =  \RZ_{\geq 0} ( 1 ,  k ) +
\RZ_{\geq 0} ( -1 , l  ) \ ,
\]
hence $\rho$ is always strongly convex by Proposition 1, and
$$
\begin{array}{cll}
{\cal C}^*_\rho & = & \RZ_{\geq 0} ( l , 1 ) +
\RZ_{\geq 0} (  - k , 1 ) \\
\triangle (L^*)_\rho & = & [ - k , l ] \ . 
\end{array}
$$
The 
bundle ${\cal O} ( \iota^* \rho )$ is the ample line 
bundle ${\cal O}_{\PZ^1} ( k + l )$. 
The map $p_{\rho}$ of (\req(rho)) is given by
\[
p_{\rho} : D  \longrightarrow \overline{L}  \ , \ \
( m_1 , m_2 ) \mapsto ( m_1 - m_2  , km_1 + lm_2  ) ,
\]
whose cokernel is generated by $q_L$ with $(k+l)q_L = 
p_{\rho} ( 1 , 1 )$. The dual map $p_{\rho}^*$ becomes:
\[
p_{\rho}^* : \overline{L}^* ( = \ZZ^2 ) \longrightarrow 
D^* ( = \ZZ^2 ) \ , \ ( b , a ) 
\mapsto ( ak + b , al - b  ) \ .
\] 
$p_{\rho}$ induces an isomorphism between the first 
quadrant cone in $D_{\RZ}$ and
the cone ${\cal C}_{\rho}$ in $\overline{L}_{\RZ}$.
Hence the space $E^*_{\rho}$ is the quotient of  
${\CZ}^2 - \{ 0 \}$ by the multiplcation  
$\omega$ (= the primitive 
$(k+l)$-th root of unity), and the bundle map induced
by $p_{\rho}$ is the projection
\[
{\CZ}^2 - \{ 0 \} \longrightarrow 
({\CZ}^2 - \{ 0 \}) / < \omega > \ = \ E^*_{\rho} \ .
\]
The basis of $\Gamma ( {\PZ}_{(\Sigma , L)}
, {\cal O }(\iota^* \rho))$ in Proposition 2 consists of
\[
( m , 1 ) \ , \ \ m \in \ZZ \cap [ -k , l ] \ ,  
\]
bijective to the following elements under $p_{\rho}^*$:
\[
( m_1 , m_2 ) \in \ZZ_{\geq0}^2 \ , \
m_1 + m_2 = k + l .
\]
This gives the usual monomial polynomial basis for 
$\Gamma ( {\PZ}^1, {\cal O }(k+l))$. As a consequence, for a
general hypersurface $X$ of ${\PZ}_{(\Sigma , L)}$ defined by
a section of ${\cal O }(\iota^* \rho)$, the equality
\be
| X | = | L^* \cap \mbox{Int} ( \triangle (L^*)_\rho ) | + 1 
\ 
\ele(|X|n=1)
holds.

(II) The weighted projective $n$-space 
$\PZ^n_{(n_i)}$ with weights 
$n_i, 1 \leq i \leq n+1, $ satisfying 
gcd $(n_j \ | \ j \neq i ) = 1$ for all $i$. 
Now $L$ is
the $n$-dimensional lattice generated by $n+1$ elements 
$\delta_j , 1 \leq j \leq n+1,$ with the only relation
\[
\sum_{j = 1}^{n+1} n_j \delta_j = 0 \ .
\]
Set
$$
\begin{array}{l}
D  \ = \ \ZZ^{n+1} \ , \\
\beta : D \longrightarrow L \ , \ \ \ (m_1, \ldots , m_{n+1})
\mapsto \sum_{j=1}^{n+1} m_j \delta_j \ . \\
\Sigma  =  \{ \sum_{j \in I}\RZ_{\geq 0} \delta_j 
\ | \ I \subset \{ 1, \ldots , n+1 \} \ , \ | I | \leq n \} \ , \\
\end{array}
$$ 
We have 
$$
\begin{array}{ccl}
{\bf n} & = & \ZZ ( n_1, \ldots , n_{n+1} ) \ , \\
\Sigma^{1} & = & \{ \delta_{j}  \}_{j=1}^{n+1}\ , \\
\PZ_{(\Sigma , L)} & = & \PZ^n_{(n_i)} \ ,
\end{array}
$$
and (\req(principlesigma)) is the natural projection:
\[
\pi : \CZ^{n+1} -\{ 0 \} \longrightarrow \PZ^n_{(n_i)} \ , 
\ ( z_1 , \ldots , z_{n+1} ) \mapsto 
[ z_1 , \ldots , z_{n+1} ] \ .
\]
Let us make the following identifications:
\[ D^* = \ZZ^{n+1} \ , \ \ \
L^* = \{ ( k_1, \ldots , k_{n+1} ) \in D^* \ | \ 
\sum_{j} n_j k_j = 0 \ \} \ .
\]  
For 
\[
\rho = ( \rho_1 , \ldots , \rho_{n+1} ) \ \in D^* - \{ 0 \} \ , \ 
\ \rho_j \in \ZZ_{\geq 0} \ ,
\]
the linear map 
\[
p_{\rho} : D \ ( \ = \ZZ^{n+1} ) \longrightarrow \overline{L} 
 \ , \ \
( m_1 , \ldots , m_{n+1} ) \mapsto 
( \sum_{j} m_j \delta_j , \ \sum_{j} m_j  \rho_j   ) 
\]
induces an isomorphism between
the first quadrant cone in $D_{\RZ}$ and the cone 
${\cal C}_{\rho}$ in $\overline{L}_{\RZ}$. Hence $\rho$ is 
strictly convex by Proposition 1.  
It is easy to see the relation 
\[
q_L = \sum_{j=1}^{n+1} r_j ( \rho_j , \delta_j  ) \ , \ \
\ \ r_j := \frac{n_j}{ \sum_{k=1}^{n+1} n_k \rho_k } \ 
\]
holds in $\overline{L}_{\QZ}$.
The affine variety $A_{\rho}$ of (\req(ASpec)) is isomorphic
to the quotient of $\CZ^{n+1}$ by the cyclic group generated
by  a diagonal element 
dia.$[ r_1, \ldots , r_{n+1} ]$. One 
can easily see that
$$
\begin{array}{cll}
{\cal C}^*_\rho & = & \sum_{j=1}^{n+1} \RZ_{\geq 0} ( 
(a_{j,1}, \ldots , a_{j , n+1} ) , 1 )  \ , 
\ a_{j, i}   = \left\{ \begin{array}{ll}
- \rho_i, & \mbox{for} \ i \neq j , \\
\frac{1}{n_j} \sum_{k \neq j} \rho_k n_k , & \mbox{for} \ i = j 
\end{array}
\right. \\ [2mm]
\triangle (L^*)_\rho & = & \mbox{the \ convex \ hull \ spanned \ 
by } \ (a_{j,1}, \ldots , a_{j , n+1} ) \ , \ 1 \leq j \leq n+1 . 
\end{array}
$$
Therefore the integral condition for $\triangle (L^*)_\rho$ 
is equivalent to 
\[
n_j | \ \sum_{k=1}^{n+1} \rho_k n_k  \ \ \mbox{for \ all \ } 
j \ . 
\]

(III) Toric variety which dominates $\PZ^n_{(n_i)} / G $
 with $G$  
a finite diagonal linear group of  
$ \CZ^{n+1}$. Consider $\CZ^{* n+1}$ as the diagonal group and
denote 
\[
\mbox{exp} : \RZ^{n+1} \longrightarrow \CZ^{* n+1} \ , \
(x_1, \ldots , x_{n+1}) \mapsto (e^{2 \pi i x_1}, \ldots , 
e^{2 \pi i x_{n+1}}) \ .
\]
$\mbox{exp}^{-1} ( G )$ is a lattice in $\RZ^{n+1}$ containing the
standard one $\ZZ^{n+1}$, and $\mbox{exp}^{-1} ( G ) \cap \  
\QZ (n_1, \ldots , n_{n+1})$ is a 1-dimensional lattice 
with the generator $(q_1, \ldots , q_{n+1})$. Let $L$ be the 
$n$-dimension lattice defined by
\[
L := \mbox{exp}^{-1} ( G ) / \ZZ (q_1, \ldots , q_{n+1})  \ ,
\]
and $\delta_j , 1 \leq j \leq n+1 ,$ be the elements in $L$
corresponding to the standard basis of 
$\ZZ^{n+1}$. Define the fan $\Sigma_0$ in $L_{\RZ}$ by
\[
\Sigma_0 = \{ \ \sum_{j \in J} \RZ_{\geq 0} \delta_j \ | \
J \subset \{ 1 , \ldots, n+1 \} , \ | J | \leq n \ \} \ .
\]
Then $\PZ^n_{(n_i)} / G $ is the toric variety 
$\PZ_{( \Sigma_0 , L )}$. A refinement $\Sigma$ of
$\Sigma_0$ gives rise a toric morphism
\[
\PZ_{( \Sigma , L )} \longrightarrow \PZ_{( \Sigma_0 , L )} 
\ .
\]
Irreducible toric divisors of $\PZ_{( \Sigma_0 , L )}$ are parametrized
by
\[
\Sigma_0^1 = \{ \delta_1 , \ldots , \delta_{n+1} \} \ ,
\]
and they give rise to $n+1$ irreducible toric divisors in 
$\PZ_{( \Sigma , L )}$, in which 
the rest ones are given by $e_{\delta} , \delta \in 
\Sigma^{1} - \Sigma_0^{1}$. For $\delta \in 
\Sigma^{1} - \Sigma_0^{1}$, one has the expression
\[
\delta = \sum_{j=1}^{n+1} a_{\delta , j} \delta_j \ , \
a_{\delta , j} \in \QZ \ . 
\]
Then the subspace ${ \bf n }_{\CZ}$ of 
$D_{\CZ}$ is generated by the following elements:
\[
\sum_{j=1}^{n+1} q_j e_{\delta_j} \ ; \ \ \ 
e_{\delta} - \sum_{j=1}^{n+1} a_{\delta , j} e_{ \delta_j } \ ,
\ \delta \in \Sigma^{1} - \Sigma_0^{1} \ , 
\]
which form a basis of ${ \bf n }_{\CZ}$ 
indexed by $\{ 0 \} \cup \Sigma^{1} - \Sigma_0^{1}$.
Hence one has the following isomorphic vector spaces: 
\be
\mbox{Pic} ( \PZ_{( \Sigma , L )} )_{\CZ} \ \ \simeq \ \ 
{ \bf n }^*_{\CZ} \ \ \simeq \ \ 
\bigoplus_{ z \in \{ 0 \} \cup \Sigma^{1} - \Sigma_0^{1} } 
\CZ z \ .
\ele(nbasis)

\section{ Rational Picard groups of hypersurfaces }
Now we are going to determine the rational Picard groups of
hypersurfaces in a toric varieties.  The non-ample 
hypersurfaces will be our main concern here. 
For the rest of this paper, we shall always assume the dimension of
$\PZ_{( \Sigma , L)}$ is at least 4, i.e.
\[ 
 n \ \geq \ 4 \ .
\]
First we consider the case when ${\cal C}_\rho$ is strongly 
convex, ( a weaker condition than the ampleness of 
${\cal O} ( \iota^* \rho )$ ). 
A similar conclusion for the second cohomology in Proposition 3 
holds for this situation, which was stated in \cite{DK} without 
the proof. Here we present a detailed argument.   
\par \vspace{0.2in} \noindent
{\bf Proposition 4. } Let $\rho$ be an element in $D^*$ 
such that ${\cal C}_{\rho}$ is a 
strongly convex cone with $\{ ( \delta , \rho_\delta ) \ | \ 
\delta \in \Sigma^{1} \}$ as the fundamental generators.
Let $X$ be 
a quasi-smooth hypersurface of $\PZ_{(\Sigma, L)}$ defined by a 
section of ${\cal O} ( \iota^* \rho )$. Then  
the following vector spaces are naturally isomorphic:
\[
H^2 ( \PZ_{( \Sigma , L)} , \CZ ) \  \simeq  \  
H^2 ( X , \CZ ) \  \simeq  \mbox{Pic} ( X )_{\CZ} \  .
\]
Consequently the inclusion map induces the isomorphism: 
\[
\mbox{Pic} \ ( \PZ_{( \Sigma , L)})_{\CZ} \ \simeq \ 
\mbox{Pic} ( X )_{\CZ} \ .
\]
\par \vspace{0.2in} \noindent
{\it Proof. } Let $U$ be the complement of $X$ in 
$\PZ_{( \Sigma , L)}$. By the argument in Proposition 3, we have 
the exact sequence (\req(MVseq)) and the isomorphism 
(\req(Poincare)). For the first isomorphism of 
the conclusion, it suffices to show the vanishing of the
following cohomology (with complex coefficients 
throughout):
\be
H^{2n-3} ( U ) = H^{2n-2} ( U ) = 0 \ .
\ele(320) 
Let $\Sigma_0 $ be the
complete fan in $L_{\RZ}$ obtained by the projection of 
proper faces of ${\cal C}_{\rho}$. Note that $\Sigma_0 $ 
may not be simplical in general, while $\Sigma$ is a simplicial refinement
$\Sigma_0 $ with $\Sigma_0^{(1)} = \Sigma^{(1)}$. Hence one has the
equivariant morphism:
\be
\phi : \PZ_{( \Sigma , L)} \longrightarrow \PZ_{( \Sigma_0 , L)} \ ,
\ele(Sigma0)
which induces the isomorphic toric divisor-groups. The graph of $f_\rho$
is linear on each polyhedral cone in $\Sigma_0 $. By the strong 
convexity of ${\cal C}_\rho$, $\rho$ determines an
ample toric divisor in
$\PZ_{( \Sigma_0 , L)}$,  whose
pull-back under the map $\phi$ equals to the bundle ${\cal O} ( \iota^* \rho )$ over
$\PZ_{( \Sigma , L)}$. As they have the same global sections,
there is an ample hypersurface $X_0$ of 
$\PZ_{( \Sigma_0 , L)}$ with $ X = \phi^{-1} ( X_0 )$. Denote
\[ 
U_0 : = \PZ_{( \Sigma_0 , L)} - X_0 ,
\]
then $U_0$ is an affine variety with $\phi^{-1} ( U_0 ) = U$ and 
$H^j ( U_0 ) = 0 $ for $ j > n$. Since 
$( \delta , \rho_\delta  ) \ , \ 
\delta \in \Sigma^{1} $, are the fundamental generators of ${\cal C}_\rho$,
the exceptional set of $\phi$ is of codimension greater than one.
There exists an analytic subspace $V_0$ of $U_0$ such that 
$V ( : = \phi^{-1} ( V_0) )$ is of codimension $\geq 2$ and
\[
U - V \stackrel{\phi}{ \simeq } U_0 - V_0 . 
\]
For a suitable "tubular" neighborhood $T$ of $V$ in $U$, and
$T^\circ = T - V$, we have the 
following commutative Mayer-Vietoris exact sequences:
$$
\begin{array}{lccccl}
 \longrightarrow & H^{2n-3} ( U_0 ) 
\longrightarrow & H^{2n-3} ( U_0 - V_0 ) \oplus H^{2n-3} ( V_0 ) 
\longrightarrow 
&H^{2n-3} ( \phi(T^\circ) ) \longrightarrow &H^{2n-2} ( U_0 ) &
\longrightarrow  \\
& \downarrow & \downarrow & \downarrow & \downarrow &\\
 \longrightarrow & H^{2n-3} ( U ) 
\longrightarrow & H^{2n-3} ( U - V ) \oplus H^{2n-3} ( V ) 
\longrightarrow 
&H^{2n-3} ( T^\circ ) \longrightarrow &H^{2n-2} ( U ) &
\longrightarrow 
\end{array}
$$
Since
\[
H^j ( U_0 ) = H^j ( V_0 ) = H^j ( V ) = 0 , \ \ \mbox{for} 
\ \ j \geq 2n-3 \ ,
\]
we obtain (\req(320)), hence the isomorphism between
$H^2 ( \PZ_{( \Sigma , L)} , \CZ )$ and $H^2 ( X , \CZ )$. 
Next we are going to show the vanishing property of 
the cohomology groups,
\[
H^1 ( X , {\cal O} ) = H^2 ( X , {\cal O} ) = 0 \ ,
\]
which implies the isomorphism between 
$H^2 ( X , \CZ )$ and $ \mbox{Pic} ( X )_{\CZ}$.
Let $L_1$ be a sublattice  of $L$ such that the
pull-back of $X$, denoted by $X_1$, is Cartier under the finite abelian cover:
\[
\PZ_{(\Sigma , L_1)} \longrightarrow \PZ_{(\Sigma , L)} \ .
\]
Now it suffices to show the vanishing of 
$H^i ( X_1 , {\cal O} )$, for $i=1,2$. By replacing $X_1, 
L_1$ by $X$ and $L$, we may 
assume the hypersurface $X$ to be Cartier in 
$\PZ_{(\Sigma , L)}$, but with no quasi-smooth condition required, 
 for the purpose of vanishing property of $
H^i ( X , {\cal O} )$ for $i =1, 2$. 
By the cohomology sequence associated to the exact sequence 
of sheaves over $\PZ_{(\Sigma , L)}$:
\[
0 \longrightarrow {\cal O}_{\PZ_{( \Sigma , L)}} ( - X )
\longrightarrow {\cal O}_{\PZ_{( \Sigma , L)}} 
\longrightarrow {\cal O}_X \longrightarrow 0 \ ,
\]
and the well-known fact
\[
H^j ( \PZ_{( \Sigma , L)} , {\cal O} ) = 0 \ \ \ \mbox{for}  
\ j > 0 \ ,
\]
the vanishing of $
H^i ( X , {\cal O} ), ( i =1, 2 )$ follows from
\[
H^j ( \PZ_{( \Sigma , L)} , {\cal O}_{\PZ_{( \Sigma , L)}} ( - X ) ) 
= 0 \ \ \ \mbox{for} \ j = 2 , 3 . \ 
\]
Note that ${\cal O}_{\PZ_{( \Sigma , L)}} ( - X )$ is the sheaf of 
sections for the line bundle ${\cal O} ( - \iota^* \rho )$.
By Serre duality, $
H^j ( \PZ_{( \Sigma , L)} , {\cal O}_{\PZ_{( \Sigma , L)}} ( - X ) ) $
is dual to $H^{n-j} ( \PZ_{( \Sigma , L)} , \  
\omega_{\PZ_{( \Sigma , L)}} ( X ) )$, where 
$\omega_{\PZ_{( \Sigma , L)}}$ is the canonical sheaf of 
$\PZ_{( \Sigma , L)}$. Under
the map $\phi$ of (\req(Sigma0)), we have
\[
\omega_{\PZ_{( \Sigma , L)}} ( X ) = 
\phi^*(\omega_{\PZ_{( \Sigma_0 , L)}} ( X_0 )) .
\]
Since $X_0$ is an ample divisor in $\PZ_{( \Sigma_0 , L)}$,
\[
H^k ( \PZ_{( \Sigma_0 , L)} , \  
\omega_{\PZ_{( \Sigma_0 , L)}} ( X_0 ) ) = 0 \ \ \ \ \mbox{for} \
k > 0 .
\]
By $n \geq 4$ and
\[
R^k_{ \phi_* }  {\cal O }_{\PZ_{( \Sigma , L)}}  = 0 \ \ \ \mbox{for} \
k > 0 \ , 
\]  
one has 
\[
H^{n-j} ( \PZ_{( \Sigma , L)} , \  
\omega_{\PZ_{( \Sigma , L)}} ( X ) ) \simeq 
H^{n-j} ( \PZ_{( \Sigma_0 , L)} , \  
\omega_{\PZ_{( \Sigma , L)}} ( X_0 ) ) = 0 \ \ \ \ \ 
\mbox{for} \ j = 2 , 3 .
\]
This completes the proof of this proposition.
$\Box$ \par \vspace{0.2in} \noindent
We now generalize Proposition 4 to a larger
class of hypersurfaces in toric varieties.
\par \vspace{0.2in} \noindent
{\bf Theorem 2. } Let $X$ be a quasi-smooth hypersurface of 
$\PZ_{( \Sigma , L)}$ defined by a section $s$ of 
${\cal O} ( \iota^* \rho )$ for a convex element $\rho$ in 
$D^*$ . Let $\Sigma ( \rho )$ be the fan in 
$\overline{L}_{\RZ}$ defined by (\req(Sigmarho)).   
Assume the Newton polygon of $s$ is
equal to $\triangle ( L^* )_\rho$, and $\Sigma ( \rho )$  
satisfies the following conditions:

(i)  
\[
\Sigma ( \rho )^{(1)} \cap 
\bigcup_{F : \ \mbox{facet\ of\ }{ \cal C}_{\rho}}
 \mbox{Int} ( F ) = \emptyset \ .
\]

(ii) $\Sigma ( \rho )$
is a refinement of some fan $\Lambda$ in $\overline{L}_{\RZ}$ 
with:
\[
|\Lambda | = \partial {\cal C}_{\rho} \ , \  
\Lambda^{(1)} = \{ 
\mbox{ 1-face \ of \ } {\cal C}_{\rho} \} \ .
\]
Then we have the isomorphism of vector spaces:
\[
\mbox{Pic} ( X  )_{\CZ} \simeq \mbox{Pic} 
( \PZ_{(\Sigma, L)} )_{\CZ} 
\oplus \bigoplus_{\tau , z } \{ \CZ z \ | \ z \in ( \mbox{Int} 
( \underline{\tau} ) \cap \triangle ( L ) \cap L ) \times 
( \mbox{Int} ( \underline{\tau^\prime} ) \cap 
\triangle ( L^*)_{\rho} \cap L^* ) \} \ ,
\]
where the index $\tau$ runs over the codimensional 2 faces of 
${\cal C}_{\rho}$ with its dual $\tau^\prime$ in
${\cal C}_{\rho}^*$, $\underline{\tau} , \underline{\tau^\prime}$
mean the projection on $L_{\RZ}$ , $L^*_{\RZ}$ respectively,
 and $\triangle ( L ), \triangle ( L^*)_{\rho}$ are defined by (\req(deftria)) 
(\req(tria*)). 
\par \vspace{0.2in} \noindent
{\it Proof.} Let $\Sigma_0$ be the complete fan in 
$L_{\RZ}$ obtained by the projection of $\Lambda$.
Then $\Sigma$ is a refinement of 
$\Sigma_0$, and one has the $\TZ ( L )$-morphism,
\[
\varphi : \PZ_{(\Sigma , L )} \longrightarrow \PZ_{(\Sigma_0 , L )} ,
\]
with the relation of toric divisors given by
\[
\Sigma^{1} \supset \Sigma_0^{1} \ .
\]
Let $\eta$ be the element in $D ( \Sigma_0 , L)^*_{\QZ}$ defined by
\[
\eta = \sum_{ \delta \in \Sigma_0^{1}} \eta_{\delta} e_{\delta} \ , \ \ \ 
\eta_{\delta} : = \rho_{\delta} \ \ \ \mbox{for} \ \ 
\delta \in \Sigma_0^{1}\ . 
\]
The cones ${\cal C}_{\rho}$ and ${\cal C}_{\eta}$ are the same
in $L_{\RZ}$, which implies 
\[
\rho = \varphi^* ( \eta ) \ 
\]
as divisors in $\PZ_{(\Sigma , L )}$.
Hence there exists a quasi-smooth hypersurface $X_0$ in 
$\PZ_{(\Sigma_0 , L )}$ with $\varphi^* ( X_0 ) = X $. The 
Newton polygon for the equation of $X_0$ equals to 
$\triangle ( L^* )_{\eta}$, which is the same as 
$\triangle ( L^* )_{\rho}$.
Since the Picard groups of $\PZ_{( \Sigma , L)}$ and $X$ are 
generated by divisors, hence one has the following identification of
vector spaces:
\bea(lcl)
\mbox{Pic}( X )_{\CZ} & = & \mbox{Pic}( X_0 )_{\CZ} \oplus \bigoplus \{ 
{\CZ} e 
\ | \ e : \mbox{irreducible \ exceptional\ divisor \ 
in \ } X \ \mbox{over}
\ X_0 \ \} \\ [2mm]
\mbox{Pic}( \PZ_{(\Sigma , L )} )_{\CZ} & = & 
\mbox{Pic}( \PZ_{(\Sigma_0 , L )} )_{\CZ}  \ \oplus \\
&&\bigoplus \{ {\CZ} E 
\ | \ E : \mbox{irreducible \ exceptional\ divisor\ in \ } 
\PZ_{(\Sigma , L )} \ \mbox{over} \ \PZ_{(\Sigma_0 , L )} \ \} \\
\elea(PicXP)
In the second equality above,  the right hand side 
is a vector space with a basis parametrized by 
$\Sigma^{1}$. By Proposition 4, we have
\[
\mbox{Pic}( X_0 )_{\CZ} \ = 
\ \mbox{Pic}( \PZ_{(\Sigma_0 , L )} )_{\CZ} \ .
\]
Let $\mbox{orb} ( \lambda )$
be $\TZ ( L )$-orbit in $\PZ_{(\Sigma_0 , L )}$
for $\lambda \in \Sigma_0$. The closure
$\overline{\mbox{orb} ( \lambda )}$ of $\mbox{orb} ( \lambda )$
is a toric variety 
having the lattice and the fan given by 
\[
L / \ZZ ( L \cap \lambda ) \ , \ \ 
\mbox{Star} (\lambda ) : =
\{ \overline{\sigma} \ | \ \sigma \in \Sigma_0 \ , \ \lambda 
\prec \sigma \ \} \ ,
\]
where $\overline{\sigma}$ denotes the image of $\sigma$ in 
$L_{\RZ} / \RZ \lambda$. Hence 
dim $\overline{ \mbox{orb} ( \lambda ) }$ + dim $( \lambda )$
= $n$. Every element $\lambda$ in $\Lambda$ is contained in
an unique face $\tau ( \lambda )$ of ${\cal C}_{\rho}$
with the same dimension. One has 
\[
{\cal C}^*_{\rho} \cap \lambda^{\bot} = 
\tau ( \lambda )^{\prime} \ \ ( \ = 
{\cal C}^*_{\rho} \cap \tau ( \lambda )^{\bot} \ ) \ .
\]
Since the Newton polygon for
$X_0$ is equal to $\triangle ( L^* )_{\rho}$, together with
the equality of (\req(|X|n=1)) for each 1-dimensional 
$\overline{\mbox{orb} ( \lambda )}$, one concludes the 
following relation holds:
\[
X_0 \bigcap \overline{\mbox{orb} ( \lambda )} = 
\left\{ \begin{array}{ll}
\emptyset , & \mbox{if} \ 
\mbox{dim} \ \overline{\mbox{orb} ( \lambda )} = 0 \ , \\
\mbox{a \ finite \ set \ with } \ | \mbox{Int} ( 
\underline{\tau ( \lambda )^\prime} ) \cap \triangle ( L^* )_{\rho}
 \cap L^*| +1 \ \mbox{elements} , & \mbox{if} \ 
\mbox{dim} \ \overline{\mbox{orb} ( \lambda )} = 1 \ , \\
\mbox{a \ connected \ set \ } , & \mbox{if} \ 
\mbox{dim} \ \overline{\mbox{orb} ( \lambda )} \geq 2 \ .
\end{array} \right.
\]
By the condition (i) in our assumption and (\req(PicXP)),
we obtain the result.
$\Box$ \par \vspace{0.2in} \noindent
\section{ Anti-canonical hypersurfaces and Calabi-Yau 3-folds}
In this section, we are going to apply our previous results 
to the discussion of an anti-canonical hypersurface 
$X$ in $\PZ_{( \Sigma , L )}$. Note that  
such $X$ has the trivial canonical 
sheaf. As the canonical sheaf of $\PZ_{( \Sigma , L )}$ is given 
by $\iota^* \kappa$ for $\kappa$ defined in (\req(defkappa)), 
 we now set the element $\rho$ of $D^*$ 
in Section 6  equal to 
$ - \kappa $ for the discussion of this section.
By the corollary of Proposition 1, the convexity of 
$- \kappa$ is equivalent to the convex property of the 
polytope $\triangle ( L )$ in $L_{\RZ}$. In 
this situation, $\triangle ( L )$ is integral 
with respect to $L$, and 
$\triangle ( L^* )_{-\kappa}$ is the dual polytope
of $\triangle ( L )$ in $L^*_{\RZ}$ defined by
\[
\triangle ( L )^* : = \{ y \in L^*_{\RZ} \ | \ 
< x , y > \ \geq -1  \ \ \ \forall \ x \in \triangle ( L ) 
\ \} \ .
\]
If the Newton polygon of a section of ${\cal O} ( - \iota^*
\kappa )$ is equal to 
$\triangle ( L^* )_{-\kappa}$, $( \triangle ( L ) , L )$ forms
a reflexive polytope in the sense of
Batyrev \cite{B}: 
 \par \vspace{.2in} \noindent
{\bf Definition 4.} Let $L$ be a lattice, and $\triangle$
a convex polytope in $L_{\RZ}$.
$( \triangle , L )$ is called a reflexive polytope if 
both $\triangle$ and $\triangle^*$ are integral.
$\Box$ \par \vspace{.2in} \noindent
Conversely, for a toric variety $\PZ_{(\Sigma, L)}$ with 
the reflexive polytope $( \triangle ( L ) , L )$,  
$- \kappa$ is an element in $D^*$. 
Consider a hypersurface $X$ of $\PZ_{(\Sigma, L)}$ with
$\omega_X = {\cal O}_X$, or equivalently to say, $X$ is 
defined by a section of ${\cal O} ( - \iota^* \kappa )$.  
As a corollary of Theorem 1 and 2, one has the following result:
\par \vspace{0.2in} \noindent
{\bf Proposition 5. } For $n \geq 4$, let $\PZ_{(\Sigma, L)}$ 
be a $n$-dimensional toric variety  with $( \triangle ( L ) , 
L)$ being a reflexive polytope. 
Let $X$ be 
a quasi-smooth hypersurface of $\PZ_{(\Sigma, L)}$ defined by
a section in ${\cal O} ( - \iota^* \kappa )$ whose Newton 
polygon equals to $\triangle ( L )^*$.
Assuming

(i) 
\[
\Sigma^{1} \cap 
\bigcup_{F  : \ \mbox{facet\ of\ } \triangle ( L ) } 
\mbox{Int} ( F ) = \emptyset \ .
\]

(ii) $\Sigma $
is a refinement of some fan $\Sigma_0$ in $L_{\RZ}$ 
with $\Sigma_0^{1} = \{ \mbox{vertex \ of} \ 
\triangle ( L ) \}$, 
\par \noindent
Then we have the isomorphic vector spaces:
\[
\mbox{Pic} ( X  )_{\CZ} \simeq {\bf n} ( \Sigma , L)^*_{\CZ} 
\oplus \bigoplus_{F , z } \{ \CZ z \ | \ z \in ( \mbox{Int} 
( F ) \cap \triangle ( L ) \cap L ) \times 
( \mbox{Int} ( F^\prime ) \cap \triangle ( L )^* \cap L^* ) \} \ ,
\]
where the index $F$ runs over the codimensional 2 faces of 
$\triangle ( L )$ with its dual face $F^\prime$ 
defined by
\[
F^\prime : = \{ y \in \triangle ( L )^* \ | \ < x , y > = -1 
\ , \  \mbox{for} \ x \in F   \} \ .
\] 
$\Box$ \par \vspace{0.2in} \noindent
For $n =4$, the hypersurface $X$ in the above
proposition is smooth under some "maximal" condition on $\Sigma^1$.
\par \vspace{0.2in} \noindent
{\bf Proposition 6. } Let $X$ be a hypersurface of 
$\PZ_{( \Sigma , L )}$ in Proposition 5 for $n = 4$, and assume
\[
\Sigma^1 = ( L \cap \partial \triangle ( L ) ) - 
\bigcup_{F : \mbox{facet\ of\ } \triangle ( L ) } 
  \mbox{Int} ( F ) \ .
\]
Then $X$ is a smooth CY space.
\par \vspace{0.2in} \noindent
{\it Proof.} As in the proof of Theorem 2, we consider 
the morphism $\varphi$ from 
$\PZ_{( \Sigma , L)}$ to $\PZ_{( \Sigma_0 , L)}$, and  
the quasi-smooth hypersurface $X_0$ of 
$\PZ_{( \Sigma_0 , L)}$ with $X = \varphi^* ( X_0 )$. 
Now $X_0$ has the trivial canonical sheaf and the
Newton polygon of its equation equals to $\triangle
( L )^*$. Given an element $p$ in $X_0 $, we are
going to show that $X$ is smooth near 
$\varphi^{-1} ( p )$. Let $\lambda$ be the element in $\Sigma_0$
such that the $\TZ ( L )$-orbit 
orb($\lambda$) in $\PZ_{( \Sigma_0 , L)}$ contains $p$. From the assumption on 
$\Sigma_0$ and $\Sigma$, it needs only to consider the case
when $m :=$ dim orb($\lambda ) = 2 , 3$. 
Let $\sigma$ be a 3-dimensional simplicial cone in
$\Sigma_0$ which contains $\lambda$ as proper face.
$\sigma$ is generated by the vertices $\delta_i, 1 \leq i \leq 4,$ of
$\triangle ( L )$, and one may assume
\[
\sigma = \sum_{j=1}^4 \RZ_{\geq 0} \delta_j \ , \ \
\lambda = \sum_{j=1}^m \RZ_{\geq 0} \delta_j \ .
\]
Associated to $\sigma$, there is the affine open subset 
$U_{\sigma}$ of $\PZ_{( \Sigma_0, L)}$ with a finite abelian 
cover (\req(psi)):
\[
\psi : \CZ^4 \longrightarrow U_{\sigma} = \CZ^4 / I_{\sigma} 
\ , \ \ 
U_{\sigma} := \mbox{ Spec } \CZ [ \breve{ \sigma } \cap L^* ]
, \ \ I_{\sigma} := L / \sum_{i=1}^4 \ZZ \delta_i \ .
\]
Let $(t_1, t_2, t_3 , t_4)$ be the coordinates of $\CZ^4$ 
corresponding to the dual basis $\{ \delta^*_i \}_{ i=1}^4$.
By (\req(eqnX)), 
$\psi^* ( X_0 \cap U_{\sigma})$ is a non-singular hypersurface 
in $\CZ^4$ defined by a $I_{\sigma}$-invariant polynormial:
\[
\eta ( t ) = \sum_{ y \in L^* \cap \triangle (L)^*} \alpha_{y} t^{k (y)} ,
\ \ \ t^{k (y)} := \prod_{i=1}^4 t_i^{k (y)_i} \ , \ \
{k (y)_i} := < \delta_i , y > + 1 .
\]
Let $y_0$ be the vertice of $\triangle (L)^*$ dual to the facet
of $\triangle ( L )$ containing all $\delta_i$'s. Since the
Newton polygon for $X_0$ is equal to $\triangle (L)^*$, $y_0$
belongs to $L^*$ and 
\[
\alpha_{y_0} t^{k (y_0)} = \alpha_{y_0} \neq 0 \ .
\]
In the case for $m=3$, by Lemma 4 and the condition (i) of 
the assumption, $ \delta_1, \delta_2, \delta_3$, 
together with some vertex $\delta_4^{\prime} $ of 
$\triangle ( L )$ not equal to $\delta_4$,
forms a 3-dimensional element in $\Sigma_0$, hence its
dual vertex $y_1$ in 
$\triangle ( L )^*$ is not equal to $y_0$. 
The same conclusion
holds also for the case $m=2$ if one starts with a suitable choice of $\delta_i,
1 \leq i \leq 4$. Therefore we have
\[
\alpha_{y_1} t^{k (y_1)} = \alpha_{y_1} t_4^{k_4} \ , \
k_4 \geq 2 , 
\]  
and the projection to the first three components 
$(t_1, t_2, t_3)$ defines the local coordinate system of 
$\psi^* ( X_0 \cap U_{\sigma})$ near $\psi^{-1}(p)$. Through
the functions  $t_1, t_2, t_3 $ and $ f ( t )$, one has the 
isomorphims:
\be
( \PZ_{( \Sigma, L_0 )} , p ) \ \simeq \
( \CZ^3 / G , 0 ) \times ( \CZ , 0 ) \ \simeq \
( X_0 , p ) \times ( \CZ , 0 ) \ ,  
\ele(PX0)
where $G$ is a finite diagonal subgroup of $SL_3 ( \CZ )$ 
 (by the trivial canonical sheaf of $X_0$). The group
of $G$ is described by  
\[
G \ \simeq \ 
( L \cap \sum_{j=1}^m \RZ \delta_j ) / \sum_{j=1}^m \ZZ \delta_j \ .
\]
Denote $<\delta_j>_{j=1}^m$ the convex set spanned by 
$\{ \delta_j \}_{j=1}^m$.
Since $G$ is a subgroup of $SL_m ( \CZ )$ for $ m = 2 , 3$,
$ L \cap \sum_{j=1}^m \RZ \delta_j$ is a sublattice of $L$
generated by $L \cap <\delta_j>_{j=1}^m$, whose
classes generate the group $G$.
The fan $\Sigma$ induces a triangulation of 
$<\delta_j>_{j=1}^m$ with $L \cap <\delta_j>_{j=1}^m$ as the
set of vertices, hence it gives a crepant toric resolution of
$\CZ^3 / G$ \cite{MOP} \cite{R89}. By (\req(PX0)), both $\PZ_{( \Sigma, L )}$ and
$X$ are non-singular near $\psi^{-1}(p)$. Therefore we obtain
the result of this proposition.
$\Box$ \par \vspace{0.2in} \noindent
{\bf Remark.}    
When the polytope is a 4-dimensional
simplex, the concept of reflexive simplex $( \triangle , L )$ is 
equivalent to  
${\PZ}_{(n_i)}^4/G $ in Sect. 5 (III) with $G \subset SL_5 ( \CZ)$
and the "Fermat-type" condition on weights $n_i$'s,
( a precise statement see Appendix (II) ). In this situation,
Proposition 5 and 6 are given in \cite{R91} \cite{R93y} 
\cite{R94}.

\section{Appendix}
 
(I) {\bf Remark on Proposition 3 . } For a hypersurface $X$ 
of $\PZ_{(\Sigma, , L )}$
in Proposition 3 , one has the isomorphic 
Hodge groups:
\[
H^q ( X , \Omega^p_X ) \simeq H^{q+1} ( \PZ_{(\Sigma, L )} , 
\Omega^{p+1}_{\PZ_{(\Sigma, L )}} ) \ , \ \
\mbox{for} \ p + q \geq n \ ,
\]
hence 
\[
H^q ( X , \Omega^p_X ) = 0 \ , \ \ \mbox{for} \ p + q \geq n \ , \
p \neq q \ .
\]
In fact, the spectral sequence
\[
E^{pq}_1 = H^q ( \PZ_{(\Sigma, L )} , 
\Omega^p_{\PZ_{(\Sigma, L )}} ( \log X ) ) \ \ \Longrightarrow 
H^{p+q} ( \PZ_{(\Sigma, L )} - X , \CZ ) \ ,
\]
degenerates at $E_1$-term \cite{Dan} \cite{Steen}.
By (\req(HU=0)) , we have
\[
H^q ( \PZ_{(\Sigma, L )} , 
\Omega^p_{\PZ_{(\Sigma, L )}} ( \log X ) ) = 0 \ ,
\ \ \ \mbox{for} \ p + q > n \ .
\]
Then the conclusion follows from 
\[
H^q ( \PZ_{(\Sigma, L )} , 
\Omega^p_{\PZ_{(\Sigma, L )}} ) = 0 \ ,
\ \ \ \mbox{for} \ p \neq q  \ ,
\]
and the cohmologoy sequence for the exact sequence of sheaves: 
\[
0 \longrightarrow 
\Omega^{p+1}_{\PZ_{(\Sigma, L )}} 
\longrightarrow  \Omega^{p+1}_{ \PZ_{(\Sigma, L )}} ( \log X ) 
\longrightarrow \Omega^p_X \longrightarrow 0  \ .
\]

(II) {\bf Reflexive simplex . } The following two sets are in one-one 
correspondence : \par \noindent
(i) The collection of all reflexive pair $(\triangle , M)$ 
with $\triangle$ a 4-dimensional simplex. \par \noindent
(ii) The collection of all $( {\PZ}_{(n_i)}^4, G ) $ with
the weights $n_i$'s and the group $G$ satisfying the conditions:
\[
d_j := \frac{d}{{\rm g.c.d}(d, n_j)}  \in {\ZZ} \ \ \ \forall \ j \ ,\ \  \ 
d := \sum_{i=1}^n n_i \ , \ \ \ \
  Q \subseteq G \subseteq SD 
\]
where $Q$ is the group generated by ${\rm dia.}[ e^{2\pi i n_1/d}, \ldots , e^{2\pi i n_5/ d}]$, and 
$SD = \{ {\rm dia.}[ t_1, \ldots , t_5 ] \in SL_5(\CZ) \ | \ t_i^{d_i} = 1
\mbox{ \ for \ all \ } i \} $.
\par \vspace{0.2in} \noindent
( For a proof of the above statement, see \cite{B} \cite{R94} )

\section{Acknowledgments} Part of this work was done while
the author was visiting Kyoto University and Tohoku  
University of Japan during the early spring of 1995. He 
wishes to thank 
Professors S. Mori, T. Oda and K. Saito for the 
kind invitations and warm hospitality.  
Especially to Professor T. Oda, the author was benefited a great deal  
from the many discussions, suggestions, and also the useful 
references related to this paper.

\end{document}